\documentclass[a4paper,11pt]{article}
\usepackage{amsmath,amssymb,amsthm,a4wide,color}
\usepackage[utf8]{inputenc}
\usepackage[T1]{fontenc}
\usepackage{authblk,mathrsfs}
\usepackage{algorithm,algpseudocode}
\usepackage{graphicx,subcaption}

\usepackage{hyperref}

\usepackage{booktabs}

\usepackage{caption}
\usepackage{subcaption}

\usepackage{pdfsync}
\newcommand{\RR}{\mathbb{R}}
\newcommand{\CC}{\mathbb{C}}
\newcommand{\mrm}{\mathrm}
\newcommand{\mL}{\mathrm{L}}
\newcommand{\mH}{\mathrm{H}}
\newcommand{\mJ}{\mathrm{J}}
\newcommand{\mV}{\mathrm{V}}
\newcommand{\mW}{\mathrm{W}}
\newcommand{\Id}{\mathrm{Id}}
\newcommand{\mD}{\mathrm{D}}
\newcommand{\mT}{\mathrm{T}}
\newcommand{\mS}{\mathrm{S}}
\newcommand{\mP}{\mathrm{P}}
\newcommand{\mR}{\mathrm{R}}
\newcommand{\mQ}{\mathrm{Q}}

\newcommand{\mbH}{\mathbb{H}}

\newcommand{\mbX}{\mathbb{X}}
\newcommand{\Vh}{\mathrm{V}_{h}}

\newcommand{\mbVh}{\mathbb{V}_{h}}

\newcommand{\bp}{\boldsymbol{p}}

\newcommand{\bd}{\boldsymbol{d}}

\newcommand{\bx}{\boldsymbol{x}}

\newcommand{\bn}{\boldsymbol{n}}
\newcommand{\bu}{\boldsymbol{u}}
\newcommand{\bv}{\boldsymbol{v}}

\newcommand{\dir}{\textsc{d}}
\newcommand{\neu}{\textsc{n}}

\newcommand{\lbr}{\lbrack}
\newcommand{\rbr}{\rbrack}

\newcommand{\ctru}{\boldsymbol{u}}
\newcommand{\ctrv}{\boldsymbol{v}}
\newcommand{\ctrw}{\boldsymbol{w}}
\newcommand{\ctrp}{\boldsymbol{p}}
\newcommand{\ctrq}{\boldsymbol{q}}
\newcommand{\ctrf}{\boldsymbol{f}}
\newcommand{\ctrg}{\boldsymbol{g}}

\newcommand{\Bh}{\mathrm{B}}
\newcommand{\Ah}{\mathrm{A}}
\newcommand{\Sh}{\mathrm{S}}

\newcommand{\Ts}{\mT_{\textsc{s}}}

\newcommand{\Range}{\mrm{Range}}

\newcommand{\calT}{\mathcal{T}}

\newcommand{\dof}{\mathrm{dof}}
\newcommand{\loc}{\mathrm{loc}}

\newtheorem{defn}{Definition}[section]
\newtheorem{lem}[defn]{Lemma}
\newtheorem{thm}[defn]{Theorem}
\newtheorem{prop}[defn]{Proposition}
\newtheorem{cor}[defn]{Corollary}
\newtheorem{example}[defn]{Example}
\newtheorem{asum}{Assumption}

\newtheorem{remark}[defn]{Remark}

\title{\textbf{Non-self adjoint impedance in}\\
  \textbf{Generalized Optimized Schwarz Methods}}
\author[1]{X.Claeys}

\affil[1]{Sorbonne Université, Université Paris-Diderot SPC, CNRS, 
  Laboratoire Jacques-Louis Lions}

\date{}

\begin{document}

\maketitle

\begin{abstract}
  We present a convergence theory for Optimized Schwarz Methods that rely on a
  non-local exchange operator and covers the case of coercive possibly
  non-self-adjoint impedance operators. This analysis also naturally deals 
  with the presence of cross-points in subdomain partitions of arbitrary
  shape. In the particular case of hermitian positive definite impedance,
  we recover the theory proposed in \lbrack Claeys \& Parolin,2021\rbrack.
\end{abstract}

\section*{Introduction}

Although domain decomposition (DD) literature now offers a wide panel
of well established techniques in the case of symetric positive
definite systems, many problems of high applicative relevance do not
fit this case. Wave propagation problems in particular, on which the
present article focuses, traditionally lead to indefinite systems and
can only be adressed by a much smaller set of DD methods.

While there exist overlapping domain decomposition approaches adapted
to the wave context (see \cite[Chap.11]{MR2104179},
\cite{zbMATH07248609,zbMATH05626634} and references therein) we shall
rather be interested in non-overlapping domain decomposition, also
called substructuring. We wish to pay particular
  attention to geometrical configurations involving cross-points due
  to their practical relevance. Cross-points are points where
  at least three subdomains are adjacent, or two subdomains
  meet at the external boundary of the computational domain.  In the
approach proposed by Després \cite{MR1291197} also known, together
with its variants, as Optimized Schwarz Method (OSM)
\cite{zbMATH05029831,zbMATH07020343}, the wave equation with outgoing
Robin boundary condition is solved in each subdomain, and neighbouring
subdomains are coupled by swapping Robin traces across each interface.
Robin traces involve an impedance factor the choice of which has a
strong impact on the speed of convergence of the global solution
algorithm. While this impedance factor was the wave number in the
thesis of Després, it can be chosen operator valued and exterior
Dirichlet-to-Neumann (DtN) maps were spotted as an ideal choice for
model configurations with no cross-point, see
\cite{nataf:hal-02194208}. Consequently, many contributions searched
for the best tuning of this impedance parameter, trying to
approximate exterior DtN maps typically by means of Pade approximants.
In this direction we refer in particular to Antoine, Geuzaine and their
collaborators \cite{Boubendir2012,ElBouajaji2015,Modave2020b}. 

In the original work of Després, a convergence proof had been provided
for general geometrical configurations, but it gave no estimate
regarding the rate of convergence and was restricted to a continuous
(i.e. non-discretized) setting.  Later Collino \& Joly
\cite{MR1764190} proposed a general theoretical framework for OSM that
remained restricted to continuous settings and required no cross-point
in the subdomain partition. The question of cross-points remained a
thorny issue that hindered the development of a general theoretical
framework for the analysis of OSM, despite a few scattered
contributions that focused on the cross-point issue
\cite{Bendali2006,MR3519297,Gander2013}.

Until recently, all variants of OSM systematically enforced a coupling
between subdomains by means of the same local exchange operator that
swapped traces from both sides of each interface.  Because this local
swapping operator is not continuous when the subdomain partition
involves cross-points, in \cite{claeys2019new} we proposed to replace
it by a non-local counterpart which paved the way to a convergence
analysis similar to \cite{MR1764190} but in general geometrical
configurations where cross-points are allowed. This idea was then
extended to discrete settings in \cite{claeys2020robust} where an
explicit estimate of the convergence rate was provided for general
geometrical partitionning, regardless of the presence of cross-points.

The theorical framework of \cite{claeys2019new,claeys2020robust} holds
for a whole family of impedance operators.  It relies on the
hypothesis that the impedance is hermitian positive
definite (HPD), see Assumption 5.1 in \cite{claeys2020robust}. This
covers certain OSM approaches that pre-existed in the
literature, including the original Després algorithm for which a novel
convergence estimate was derived, see Example 11.4 in
\cite{claeys2020robust}.  However many variants of OSM involve
impedance operators that are not HPD. Here are a few examples:
\begin{itemize}
\item[$\bullet$] evanescent mode damping algorithm (EMDA)
  \cite[Eq.(32)]{MR2396903}, \cite[Eq.(9)]{Boubendir2012},\\[-20pt]
\item[$\bullet$] optimized Robin conditions (OO0)
  \cite[Section 3.1]{MR1924414},\\[-20pt]
\item[$\bullet$] optimized 2nd order conditions (OO2)
  \cite[Section 3.2]{MR1924414},\\[-20pt]
\item[$\bullet$] "two-sided" optimized conditions, see
  \cite[Section 3]{MR2344706},\\[-20pt]
\item[$\bullet$] "square-root" conditions
  \cite[Eq.(10)]{Boubendir2012} and its Pade approximations
  \cite[Section 5]{Boubendir2012}.
\end{itemize}
Other examples of non-self adjoint impedance operators can be found in
\cite[\S 5.2]{LECOUVEZ2014403}, \cite[\S 2.1.2 and Eq. (3.1)]{zbMATH07197799},
\cite{despres:hal-02612368}, \cite[Section 2.2]{despres:hal-03230250}.
In spite of their analogies with Després' initial algorithm \cite{MR1291197}, the strategies
listed above cannot be analyzed through the theory of \cite{claeys2019new,claeys2020robust}
because these impedance operators are not HPD. Exterior Dirichlet-to-Neumann maps
that have emerged as an ideal choice of impedance do not themselves induce 
HPD impedance operators. This is our motivation for seeking to extend our framework
beyond the case of HPD impedances.

In the present contribution, we extend this theory to the case of
impedance operators that only need to be coercive and can thus admit a
non-HPD part. In this context, the non-local exchange
operator that maintains a coupling between subdomains is not
orthogonal anymore.  Besides this extension, the present contribution
also contains several other novelties.
\begin{itemize}
\item[$\bullet$] We identify the spectral bounds of the impedance
  operator that played a pivotal role in the
  convergence estimate of \cite{claeys2020robust} as the continuity
  modulus and inf-sup constant of the trace operator, see Section
  \ref{SpectralBoundsTraceOperator}.
\item[$\bullet$] In Theorem \ref{FactorizationForm} we exhibit a
  factorized form for the inverse of the skeleton formulation
  involving an oblique projector onto the space of Cauchy traces
  parallel to the space of traces that comply with transmission
  conditions, see Proposition \ref{EstimationProjectorCauchy}.
\item[$\bullet$] We show that the local swapping operator considered
  elsewhere in the literature on OSM is simply the exchange operator
  associated to the identity matrix as impedance, and we study which
  other impedances lead to this same local exchange operator, see
  Section \ref{PurelyLocalImpedance}.
\item[$\bullet$] We show that the OSM strategies presented in
  \cite{zbMATH07197799,despres:hal-03230250} are particular cases of
  the non-HPD theory presented here.
\end{itemize}
The present contribution aims at extending the theory of
\cite{claeys2019new,claeys2020robust}.  To our knowledge this is the
first contribution that establishes a convergence theory for Optimized
Schwarz Methods for such a general class of impedance operators and
geometric partitions including the possibility of cross-points.
In the same direction, we should also point to the preprint
  \cite{arxiv.2204.03436} posted during the revision phase of the present contribution
that elaborates a general framework for substructuring methods that
covers the possibility of cross-points as well. Of course, in the special case
of symetric positive definite impedance, we recover the framework of
\cite{claeys2020robust}. We believe that
this unifying framework may help better understanding the convergence
properties of the substructuring strategies belonging to the OSM
family. We provide a numerical illustration in the
last section, although our goal is mainly theoretical. For more extensive
numerical results in the HPD case, we refer the reader
to Section 14 of \cite{claeys2020robust} and to Chapter 11 of \cite{parolin:tel-03118712}
in 2D and 3D for acoustics and electromagnetics in both homogeneous and heterogeneous media.

\section{Problem under study}\label{PbUnderStudy}
In the present section we describe the problem to be solved, as well
as basic notations regarding the discretization strategy. To present
the theoretical novelties, we base our presentation on a problem and a
geometric configuration very close to the ones of
\cite{claeys2020robust}.

\subsection{Wave propagation problem}
We consider a boundary value problem modeling scalar wave propagation
in an \textit{a priori} heterogeneous medium in $\RR^{d}$ with $d=1,2$
or $3$. The computational domain $\Omega\subset \RR^{d}$ is assumed
bounded and polygonal ($d=2$) or polyhedral ($d=3$). The material
characteristics of the propagation medium will be represented by two
measurable functions satisfying the following.

\begin{asum}\label{Hypo1}\quad\\
  The functions $\kappa: \Omega\to \CC$ and $\mu:\Omega\to (0,+\infty)$ satisfy 
  \begin{itemize}
  \item[(i)]  $\mrm{\Im m}\{\kappa(\bx)\}\geq 0$, $\mrm{\Re e}\{\kappa(\bx)\}\geq 0\;
    \forall \bx\in \Omega$\\[-20pt]
  \item[(ii)]  $\kappa_*:=\max(1,\sup_{\bx\in\Omega}\vert \kappa(\bx)\vert) <+\infty$\\[-20pt]
  \item[(iii)] $\sup_{\bx\in\Omega}\vert \mu(\bx)\vert + \vert \mu^{-1}(\bx)\vert<+\infty$.
  \end{itemize}
\end{asum}

\noindent 
These are both general and physically reasonable
assumptions. Condition \textit{(i)} above implies in particular that
$\mrm{\Im m}\{\kappa^{2}(\bx)\}\geq 0$ for all $\bx\in \Omega$. It
means that the medium can only absorb or propagate energy. In
addition, we consider source terms $f\in \mL^{2}(\Omega)$ and $g\in
\mL^{2}(\partial\Omega)$. The boundary value problem under
consideration will be
\begin{equation}\label{InitialPb}
  \left\{\;\begin{aligned}
    & \text{Find\;$u\in \mH^{1}(\Omega)$ such that}\\
    & \mrm{div}(\mu\nabla u) + \kappa^{2}u = -f\quad \text{in}\;\Omega,\\
    & \mu\partial_{\bn}u = g\quad \text{on}\;\partial\Omega. 
  \end{aligned}\right.
\end{equation}
where $\bn$ refers to the outward unit normal vector to
$\partial\Omega$.  As usual, for any domain $\omega\subset \RR^{d}$,
the space $\mH^{1}(\omega)$ refers to those $v\in \mL^{2}(\omega)$
such that $\nabla v\in\mL^{2}(\omega)$. It will be equipped with the
frequency dependent norm
\begin{equation*}
 \Vert v\Vert_{\mH^{1}(\omega)}^{2}:=\Vert \nabla  v\Vert_{\mL^{2}(\omega)}^{2} +
    \kappa_*^2\,\Vert v \Vert_{\mL^{2}(\omega)}^{2}.
\end{equation*}
As usual, Problem~\eqref{InitialPb} can be put in variational form:
Find\;$u\in \mH^{1}(\Omega)$ such that $a(u,v) = \ell(v)$ $\forall
v\in \mH^{1}(\Omega)$ where
\begin{equation}\label{BilinearFormPb}
  \begin{array}{rl}
      a(u,v):= & \hspace{-0.25cm} \int_{\Omega}\mu \nabla u\cdot
    \nabla \overline{v} -\kappa^{2} u\overline{v} \,d\bx\\[3pt]
        \ell(v):= & \hspace{-0.25cm} \int_{\Omega}f\overline{v} \;d\bx +
        \int_{\partial\Omega} g\overline{v}\;d\sigma
  \end{array}  
\end{equation}

\subsection{Discrete formulation}
We are interested in the numerical solution to this problem by means
of a classical nodal finite element scheme. We consider a regular
triangulation $\mathcal{T}_{h}(\Omega)$ of the computational domain
$\overline{\Omega} = \cup_{\tau\in\mathcal{T}_{h}(\Omega)}\overline{\tau}$.
Shape regularity of this mesh is \textit{not} needed for the subsequent
analysis. We denote $\Vh(\Omega)\subset \mH^1(\Omega)$ a space of
$\mathbb{P}_{k}$-Lagrange finite element functions constructed on
$\mathcal{T}_{h}(\Omega)$,
\begin{equation*}
  \begin{aligned}
    & \Vh(\Omega) := \{ v\in \mathscr{C}^0(\overline{\Omega}):\;v\vert_{\tau}\in \mathbb{P}_k(\overline{\tau})
    \;\forall \tau\in \mathcal{T}_{h}(\Omega)\}\\
    & \text{where}\quad \mathbb{P}_k(\overline{\tau}):=\{\text{polynomials of order $\leq k$ on $\overline{\tau}$} \}
  \end{aligned}
\end{equation*}
If $\omega\subset \Omega$ is any open subset that is resolved by the triangulation
i.e. $\overline{\omega} = \cup_{\tau\in \mathcal{T}_{h}(\omega)}\overline{\tau}$,
where $\mathcal{T}_{h}(\omega)\subset \mathcal{T}_{h}(\Omega)$, then
we denote $\Vh(\omega):=\{ \varphi\vert_{\omega},\;\varphi\in \Vh(\Omega)\}$
and also consider finite element spaces on boundaries
$\Vh(\partial\omega):=\{ \varphi\vert_{\partial\omega},\;\varphi\in \Vh(\Omega)\}$.
We will focus on the discrete variational formulation
\begin{equation}\label{DiscreteVF}
  \begin{array}{l}
    \text{Find}\;u_h \in\Vh(\Omega)\;\text{such that}\\
    a(u_h,v_h) = \ell(v_h)\quad \forall v_h\in\Vh(\Omega).
  \end{array}
\end{equation}
Devising domain decomposition algorithms to solve this discrete problem
is the main goal of the present article. This is why we need to assume that it
admits a unique solution, which is equivalent to assuming that the corresponding
finite element matrix is invertible.
\begin{asum}\label{Hypo2}\quad\\
\begin{equation}\label{UniformDiscreteInfSup}
  \alpha_{h}:=
  \mathop{\textcolor{white}{p}\inf}_{u\in \Vh(\Omega)\setminus\{0\}}
  \sup_{v\in \Vh(\Omega)\setminus\{0\}}\frac{\vert a(u,v)\vert}
    {\Vert u\Vert_{\mH^{1}(\Omega)}\Vert v\Vert_{\mH^{1}(\Omega)}}>0.
\end{equation}
\end{asum}
\quad\\
Although in practice this constant is uniformly bounded from below
$\liminf_{h\to 0}\alpha_h>0$, this uniform bound is not required for our analysis.

\subsection{Notation conventions}

Here we fix a few notation conventions that will be used all through
the manuscript.  Since the analysis presented here is fully discrete,
we will only deal with finite dimensional function spaces.
The vector spaces we shall consider will
systematically be assumed complex valued i.e. $\CC$ will always be
the scalar field.

If $\mV$ is any finite dimensional vector space
equipped with the norm $\Vert \cdot\Vert_{\mV}$, then $\mV^*$ will
refer to its topological dual i.e. the space of linear functionals
$p:\mV\to \CC$ continuous in the norm $\Vert \cdot\Vert_{\mV}$,
equipped with the norm
\begin{equation}\label{DualNormDefinition}
  \Vert p\Vert_{\mV^*} = \sup_{v\in \mV\setminus\{0\}}\vert \langle
  p,v\rangle\vert/\Vert v\Vert_{\mV},
\end{equation}
the duality pairing between two dual spaces being systematically
denoted $\langle\cdot,\cdot\rangle$. When such pairing brackets are
used, it shall be clear from the context which duality pair of spaces
$(\mV,\mV^*)$ is considered. For $v\in\mV$ and $p\in \mV^*$ we shall
write indistinctly $\langle p,v\rangle = \langle v,p\rangle$.  We
emphasize that duality pairings involve no conjugation operation.  The
polar set of any subspace $\mW\subset\mV$ will be defined by
\begin{equation*}
  \mW^\circ:=\{ \varphi\in \mV^*:\;\langle \varphi, v\rangle =
  0\;\forall v\in \mW\}.
\end{equation*}
This polar space inherits the norm \eqref{DualNormDefinition} from
$\mV^*$.  In addition if $\mL:\mV_1\to \mV_2$ is any continuous linear
map between two vector spaces $\mV_1,\mV_2$, its adjoint will then be
a map $\mL^*:\mV_2^*\to \mV_1^*$ defined by the identity
\begin{equation*}
  \langle \mL(u),\overline{v}\rangle = \langle
  u,\overline{\mL^{*}(v)}\rangle \quad \forall u\in\mV_1, \;\;\forall v\in
  \mV_2^{*}.
\end{equation*}
We shall systematically consider finite dimensional vector spaces
stemming from a finite element discretization of function spaces. The
linear operators $\mL:\mV_1\to \mV_2$ will simply be finite element
matrices but, very much like \cite[Chap.5]{zbMATH06154429}, we will
refer to them as "operators" so as to keep track of our abstract
setting where the choice of appropriate norms matters.\\
\begin{figure}[!]
  \begin{center}
    \begin{subfigure}[b]{0.25\textwidth}
      \includegraphics[height=4cm]{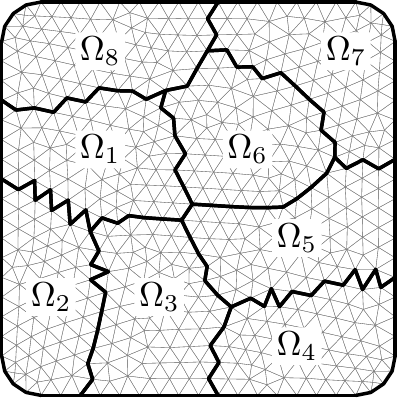}
      \caption{}
    \end{subfigure}
    \hspace{2cm}
    \begin{subfigure}[b]{0.25\textwidth} 
      \includegraphics[height=4cm]{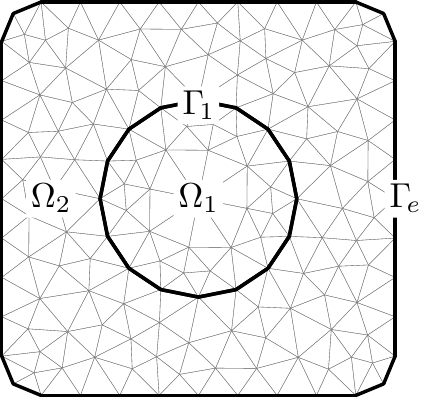}
      \caption{}
    \end{subfigure}
    \caption{Two non-overlapping decompositions of the same domain
      (a) with cross-points (b) without cross-point.}\label{MainFig}
  \end{center}
\end{figure}

\section{Geometric partitioning}\label{GeomPart}

We shall study a substructuring domain decomposition strategy for the solution
to Problem~\eqref{DiscreteVF}, which leads to introducing a non-overlapping
subdomain partition of the computational domain.  
\begin{equation}\label{SubdomainPartition}
  \begin{array}{ll}
    \overline{\Omega} = \cup_{j=1}^{\mJ}\overline{\Omega}_{j},
    & \text{with}\quad\Omega_{j}\cap\Omega_{k} = \emptyset\quad
      \text{for}\;j\neq k\\[5pt]
    \Sigma := \cup_{j=1}^{\mJ}\Gamma_{j},
    & \text{where}\quad
        \Gamma_{j}:=\partial\Omega_{j},
  \end{array}
\end{equation}
where each $\Omega_{j}\subset \Omega$ is itself a polyhedral domain
that is exactly resolved by the triangulation. We do not make any
further assumption regarding the subdomain partitionning.  As opposed
to \cite[\S 2.5.2]{MR3013465} or \cite[\S 4.2]{MR2104179}, we do not
formulate any regularity assumption on the subdomains.
Cross-points are points where at least three
subdomains are adjacent, or two subdomains meet at the external
boundary of the computational domain. For problems posed in 3D, such
points form the so-called wire basket \cite[\S
  4.6]{MR2104179}.  Our geometrical framework covers in particular
the possibility of cross-points. The treatment of such points has been
the subject of several recent contributions 
\cite{Modave2020b,Gander2013,despres:hal-03230250,claeys2020robust,claeys2019new}.
In accordance with the notations of the previous section, we have
\begin{equation*}
  \Vh(\Sigma):=\{v\vert_{\Sigma}: v\in\Vh(\Omega)\}
\end{equation*}
which is here a space of (single valued) finite element functions
defined over the skeleton that is a surface with multiple branches
i.e. the union of all interfaces.  Next we introduce continuous and
discrete function spaces naturally associated to the multi-domain
setting
\begin{equation}\label{ContinuousNorms}
  \begin{aligned}
    & \mbH^1(\Omega):= \mH^1(\Omega_1)\times \dots\times \mH^1(\Omega_\mJ),\\
    & \mbVh(\Omega) := \Vh(\Omega_1)\times \dots\times \Vh(\Omega_\mJ)\subset \mbH^1(\Omega),\\
    & \mbX_h(\Omega):= \{ (v\vert_{\Omega_1},\dots, v\vert_{\Omega_\mJ})\in\mbVh(\Omega):\;  v\in \Vh(\Omega)\}.\\[5pt]
  \end{aligned}
\end{equation}
Since they are cartesian products, these spaces are made of tuples of (volume based) functions.
The "broken space" $\mbVh(\Omega)$ is naturally identified with those functions that are piecewise
$\mathbb{P}_k$-Lagrange in each subdomain whereas, due to the matching conditions, the space 
$\mbX_h(\Omega)$ is naturally identified with $\Vh(\Omega)$ i.e. those functions that are
globally $\mathbb{P}_k$-Lagrange in the whole computationnal domain, including through interfaces
$\Gamma_j\cap \Gamma_k$. These spaces will be equipped with the norm
\begin{equation}\label{VolumeH1Norm}
  \Vert u\Vert_{\mbH^{1}(\Omega)}^{2} := \Vert u_1\Vert_{\mH^{1}(\Omega_{1})}^{2}+
  \dots+ \Vert u_\mJ \Vert_{\mH^{1}(\Omega_{\mJ})}^{2}
\end{equation}
for $u = (u_1,\dots,u_\mJ)$. Since we are interested in domain decomposition where behaviour of
functions at interfaces play a crucial role, we need to consider spaces consisting in tuples of
trace functions. These will be called Dirichlet multi-trace (resp. Dirichlet single-trace)
space 
\begin{equation}\label{MTFSpace}
  \begin{aligned}
    & \mbVh(\Sigma):= \Vh(\Gamma_{1})\times\cdots\times \Vh(\Gamma_{\mJ})\\
    & \mbX_h(\Sigma):= \{ (p\vert_{\Gamma_1},\dots,p\vert_{\Gamma_\mJ})\in\mbVh(\Sigma):\; p\in\Vh(\Sigma)\}.
  \end{aligned}
\end{equation}
The elements of $\mbX_h(\Sigma)$ can also be characterized as those tuples $\ctrp = (p_1,\dots,p_\mJ)
\in\mbVh(\Sigma)$ that match across interfaces $p_j = p_k$ on $\Gamma_j\cap\Gamma_k$.
The spaces \eqref{MTFSpace} can be obtained by taking interior traces of functions belonging
to $\mbVh(\Omega)$ resp. $\mbX_h(\Omega)$. This motivates the introduction of
the trace map $\Bh$ defined as follows
\begin{equation}\label{TraceOperator}
  \begin{aligned}
    & \Bh:\mbVh(\Omega)\to \mbVh(\Sigma)\\
    & \Bh(v):= (v_1\vert_{\Gamma_{1}},\dots,v_\mJ\vert_{\Gamma_{\mJ}})
  \end{aligned}
\end{equation}
for $v = (v_1,\dots,v_\mJ)\in\mbVh(\Omega)$.  This trace operator \eqref{TraceOperator} surjectively maps $\mbVh(\Omega)$ onto $\mbVh(\Sigma)$,
and it is also surjective from $\mbX_h(\Omega)$ onto $\mbX_h(\Sigma)$. We emphasize that the boundary trace
map \eqref{TraceOperator} is  subdomain-wise block-diagonal. Since we are in a finite dimensional
context, and $\overline{\Bh(v)} = \Bh(\overline{v})\;\forall v\in \mbVh(\Omega)$, according
to \cite[Thm. 4.7 \& 4.12]{zbMATH01022519} we have $\mrm{Range}(\Bh^*) = \mrm{Ker}(\Bh)^\circ :=
\{\phi\in \mbVh(\Omega)^*:\; \langle \phi,v\rangle = 0\;\forall v\in \mrm{Ker}(\Bh)\}$, which can be
rephrased in more concrete terms as follows. 

\begin{lem}\label{LiftingBoundaryFunctional}\quad\\
  Consider $\phi\in \mbVh(\Omega)^*$ satisfying $\langle \phi,u\rangle = 0$ for all $u\in \mbVh(\Omega)$
  with $\Bh(u) = 0$. Then there exists $\ctrp\in \mbVh(\Sigma)^*$ satisfying $\langle \phi,v\rangle =
  \langle \ctrp,\Bh(v)\rangle\;\forall v\in \mbVh(\Omega)$.
\end{lem}

\noindent 
Observe that $\mbX_h(\Omega) = \{ u\in \mbVh(\Omega):\;\Bh(u)\in \mbX_h(\Sigma)\}$.
As a consequence, if $u\in \mbVh(\Omega)$ and $\Bh(u) = 0$ then $u\in \mbX_h(\Omega)$. In other
words $\mrm{Ker}(\Bh) \subset \mbX_h(\Omega)$. The single-trace space $\mbX_h(\Sigma)$ can be
parameterized by means of the restriction operator
\begin{equation}
  \begin{aligned}
    & \mrm{R}:\Vh(\Sigma)\to \mbX_h(\Sigma)\subset \mbVh(\Sigma)\\
    & \mrm{R}(p):=(p\vert_{\Gamma_1},\dots,p\vert_{\Gamma_\mJ})\quad \text{where}\\
    & \Vh(\Sigma) :=\{ v\vert_{\Sigma}: v\in \Vh(\Omega)\} 
  \end{aligned}
\end{equation}
The space $\mbVh(\Sigma)^*$ shall be referred to as the Neumann multi-trace space,
while the following polar set will be called Neumann single-trace space
\begin{equation}\label{PolarSpace}
  \mbX_h(\Sigma)^\circ:=\{\ctrp\in \mbVh(\Sigma)^*:
  \langle \ctrp,\ctrv\rangle = 0\;\;\forall\ctrv\in \mbX_h(\Sigma)\}.
\end{equation}
This later space yields a variational characterization of $\mbX_h(\Sigma)$
through a polarity identity i.e. for any $\ctrv\in \mbVh(\Sigma)$ we have
$\ctrv\in\mbX_h(\Sigma)\iff   \langle \ctrp,\ctrv\rangle = 0\;
\forall\ctrp\in \mbX_h(\Sigma)^\circ$, which is the discrete
counterpart of a polarity result frequently used in Multi-Trace theory
\cite[Prop.2.1]{MR3069956}.  Since the range of the
operator $\mR$ is $\mbX_h(\Sigma)$ and $\overline{\mR(p)} = \mR(\overline{p})\;
  \forall v\in \Vh(\Sigma)$, this polarity identity also rewrites
\begin{equation}\label{Polarity}
  \begin{aligned}
    \mbX_h(\Sigma)^{\phantom{\circ}}      & = \mrm{Range}(\mR),\\
    \mbX_h(\Sigma)^\circ & = \mrm{Ker}(\mR^*).
  \end{aligned}
\end{equation}

\section{Reformulation of transmission conditions}\label{ReformulationTransmissionConditions}

Transmission conditions are a crucial ingredient of any domain decomposition
strategy. As a consequence, we pay a special attention to the matching conditions
at interfaces between subdomains. We need to introduce a so-called "impedance" operator $\mT$
(following the terminology of \cite{Boubendir2012}) that shall play a central role in the
subsequent analysis. All that needs to be assumed concerning this operator is the following.

\begin{asum}\label{CoercivityImpedance}\quad\\
  The linear operator $\mT:\mbVh(\Sigma)\to \mbVh(\Sigma)^*$ satisfies
  $\Re e\{\langle \mT(\ctrp),\overline{\ctrp}\rangle\}>0\;\forall \ctrp\in \mbVh(\Sigma)\setminus\{0\}$.
\end{asum}

\noindent 
We underline that, concerning the impedance operator $\mT$, no other property than
Assumption \ref{CoercivityImpedance} will be needed in the subsequent analysis. This
assumption covers the possibility that $\mT\neq \mT^*$, which is new
compared to \cite{claeys2020robust}.

In the context of waves, many domain decomposition strategies belonging to the family of
Optimized Schwarz Methods have considered particular instances of non-self-adjoint
impedances \cite{MR2396903,Boubendir2012,MR1924414,MR2344706,LECOUVEZ2014403,
    zbMATH07197799,despres:hal-02612368,despres:hal-03230250}.
With Assumption \ref{CoercivityImpedance}, the analysis of the present article
can now cope with impedance conditions that remained beyond the scope of
\cite{claeys2020robust} because they relied on non-HPD operators.
Below are a few examples of such conditions.

\begin{example}[OO0, EMDA, two sided condition]\label{ExampleEMDAOO0}
Assume a decomposition in two subdomains $\overline{\Omega} =
  \overline{\Omega}_1\cup\overline{\Omega}_2$ with no cross-point (see Fig.\ref{MainFig}(b) for an example) 
  and $\kappa\in (0,+\infty)$ and $\mu=1$. Optimized Robin conditions (OO0) \cite{MR1924414},
  two sided conditions \cite{MR2344706} and EMDA \cite{MR2396903} correspond to an impedance operator of the form
  $\langle \mT(\bu),\bv\rangle = \langle \mT_1(u_1),v_1\rangle + \langle \mT_2(u_2),v_2\rangle$
  for $\bu = (u_1,u_2), v = (v_1,v_2)$ with $u_j,v_j\in\Vh(\Gamma_j)$, where
  \begin{equation*}
    \langle \mT_j(u_j),v_j\rangle = z_j\int_{\Gamma_j}u_j(\bx)v_j(\bx) d\sigma(\bx)
  \end{equation*}
  and $z_j\in\CC$ are complex numbers chosen by means of an optimization procedure.
  In the case of OO0 and EMDA we have $z_1 = z_2$, and this constraint is relaxed
  with the two-sided condition. With EMDA one has $\Re e\{z_j\} = \kappa$.
  With OO0, EMDA and two-sided impedance, we have $\Re e\{z_j\}>0$ (see \cite[Lemma 3.3]{MR1924414})
  and a priori  $\Im m\{z_j\}\neq 0$ so that Assumption \ref{CoercivityImpedance} is
  satisfied while the impedance operators are not HPD.
\end{example}

\begin{example}\label{LienLecouvezCollinoJoly}
  Under the simplifying assumptions of the previous example,
    suppose in addition that $\Gamma_2 = \Gamma_1\cup \Gamma_{e}$ where $\Gamma_{e}
  = \Gamma_2\setminus \Gamma_1 = \partial\Omega$ is the external
  boundary of the computational domain, see Fig \ref{MainFig} (b).  For $\bu = (u_1,u_2), \bv =
  (v_1,v_2)$ with $u_j,v_j\in\Vh(\Gamma_j)$, the strategy presented in
  \cite[\S 2.1.2]{zbMATH07197799} relies on impedance operators of the
  form\footnote{The treatment of the external boundary is different in
  \cite{zbMATH07197799} compared to
  \eqref{ImpedanceCollinoJolyLecouvez}.  The main point here concerns
  the treatment of the common interface $\Gamma_1 = \Gamma_1\cap
  \Gamma_2$ though. }
  \begin{equation}\label{ImpedanceCollinoJolyLecouvez}
    \langle \mT(\bu),\bv\rangle = \langle \mT_1(u_1),v_1\rangle +
    \langle \mT_1^*(u_2\vert_{\Gamma_1}),v_2\vert_{\Gamma_1}\rangle +
    \langle \mT_e(u_2\vert_{\Gamma_e}),v_2\vert_{\Gamma_e}\rangle
  \end{equation}
  where $\mT_e = \mT_e^*$ is hermitian positive definite (HPD) and
  $\mT_1+\mT_1^*$ is positive definite. This fits Assumption
  \ref{CoercivityImpedance} above.  In Section 3 of
  \cite{zbMATH07197799}, the analysis is particularized to the case
  where $\mT_1 = (1+i\gamma)\mT_{\textsc{r}}$ where $\mT_{\textsc{r}}$
  is HPD and $\gamma\in\RR$.
\end{example}

\begin{example}\label{ExampleImpedancezTR}
  As a further instructive example, we consider an operator
  $\mT_{\textsc{r}}:\mbVh(\Sigma)\to \mbVh(\Sigma)^*$ that is
  hermitian positive definite (HPD) and choose the impedance operator
  as $\mT = z\mT_{\textsc{r}}$ where $z\in \CC$. In this situation,
  $\mT$ complies with Assumption \ref{CoercivityImpedance} provided
  that $\Re e\{z\}>0$.

 An instance of such a situation is provided by  EMDA and OO0,
  taking a reference operator $\mT_{\textsc{r}}$ stemming from surface mass matrices.
  Other examples can be constructed defining, as in \cite[Sect.3]{zbMATH07197799},
  the reference operator $\mT_{\textsc{r}}$ by means of Gagliardo semi-norms or single
  layer potentials.
\end{example}

\noindent 
Assumption \ref{CoercivityImpedance} can be rephrased by stating that the symmetric part $\Ts:=(\mT+\mT^*)/2$
induces a scalar product over $\mbVh(\Sigma)$. As a consequence the operator
$\Ts^{-1}=2(\mT+\mT^*)^{-1}$ induces a scalar product over the dual space $\mbVh(\Sigma)^*$.
We define
\begin{equation} \label{DefScalarProd}
  \begin{aligned}
    & \Vert \ctrv\Vert_{\Ts^{\phantom{-1}}}^2:= \langle \Ts(\ctrv),\overline{\ctrv}\rangle\\
    & \Vert \ctrp\Vert_{\Ts^{-1}}^2:= \langle \Ts^{-1}(\ctrp),\overline{\ctrp}\rangle\\
    & \text{where}\;\;\Ts:= (\mT+\mT^*)/2.
  \end{aligned}
\end{equation}
With these definitions, for any $\ctrv\in \mbVh(\Sigma), \ctrp\in \mbVh(\Sigma)^*$, 
we have $\vert \langle \ctrv,\ctrp\rangle\vert\leq \Vert \ctrv\Vert_{\Ts}\Vert \ctrp\Vert_{\Ts^{-1}}$
and $\Vert \Ts(\ctrv)\Vert_{\Ts^{-1}} = \Vert \ctrv\Vert_{\Ts}$ and $\Vert \Ts^{-1}(\ctrp)\Vert_{\Ts} = \Vert \ctrp\Vert_{\Ts^{-1}}$.
The theory that we present here stems from a
new treatment of transmission conditions that relies on a proper characterization of Dirichlet (resp. Neumann)
single-trace space $\mbX_h(\Sigma)$ (resp. $\mbX_h(\Sigma)^\circ$). We need a first lemma.

\begin{lem}\quad\\
  For any impedance operator $\mT:\mbVh(\Sigma)\to \mbVh(\Sigma)^*$ satisfying Assumption \ref{CoercivityImpedance},
  we have the following direct sums
  \begin{equation*}
    \begin{array}{ll}
      \mbVh(\Sigma)^{\phantom{*}}  = \mbX_h(\Sigma)\oplus \mT^{-1}(\mbX_h(\Sigma)^\circ),\\[4pt]
      \mbVh(\Sigma)^* = \mbX_h(\Sigma)^\circ\oplus \mT(\mbX_h(\Sigma)).
    \end{array}
  \end{equation*}
  Moreover, if $\mT = \mT^* = \Ts$ then these direct sums are respectively $\Ts$-orthogonal and
  $\Ts^{-1}$-orthogonal.
\end{lem}
\noindent \textbf{Proof:}

We prove the result for the second direct sum only, since the proof runs completely
parallel for the first direct sum. First pick $\ctrp\in \mbX_h(\Sigma)^\circ\cap \mT(\mbX_h(\Sigma))$
so that there exists $\ctru\in \mbX_h(\Sigma)$ with $\ctrp = \mT(\ctru)$. Then we have
\begin{equation*}
  \begin{aligned}
    0
    & = 2\Re e\{\langle \ctrp,\overline{\ctru}\rangle\} = 2\Re e\{\langle \mT(\ctru),\overline{\ctru}\rangle\}
    = \langle \mT(\ctru),\overline{\ctru}\rangle + \langle \ctru,\overline{\mT(\ctru)}\rangle\\
    & = \langle (\mT+\mT^*)\ctru,\overline{\ctru}\rangle = 2\Vert \ctru\Vert_{\Ts}^2\quad \Longrightarrow\quad \ctru = 0.
  \end{aligned}
\end{equation*}
From this we conclude that $\ctrp = \mT(\ctru) = 0$ hence $\mbX_h(\Sigma)^\circ\cap \Ts(\mbX_h(\Sigma)) = \{0\}$ since
$\ctrp$ was chosen arbitrarily.

\quad\\
Next take an arbitrary $\ctrq\in \mbVh(\Sigma)^*$ and define $\ctrv$ as the only element of $\mbX_h(\Sigma)$ satisfying
$\langle \mT(\ctrv),\ctrw\rangle = \langle \ctrq,\ctrw\rangle\;\forall \ctrw\in \mbX_h(\Sigma)$. This variational problem
admits  a unique solution due to the coercivity of $\mT$ from Assumption \ref{CoercivityImpedance}.
By construction we have $\langle \ctrq - \mT(\ctrv),\ctrw\rangle = 0\;\forall \ctrw\in \mbX_h(\Sigma)$
which means that $\ctrq' = \ctrq - \mT(\ctrv)\in \mbX_h(\Sigma)^\circ$. We have just established that
$\ctrq = \ctrq'+\mT(\ctrv)\in \mbX_h(\Sigma)^\circ+\mT(\mbX_h(\Sigma))$, hence finally $\mbVh(\Sigma)^{*} =
\mbX_h(\Sigma)^\circ\oplus\mT(\mbX_h(\Sigma))$. 

\quad\\
In the case where $\mT = \mT^* = \Ts$, take arbitrary
$\ctrp\in \mbX_h(\Sigma)^\circ,\ctrq \in \Ts(\mbX_h(\Sigma))$.
We have $\ctrq = \Ts(\ctrv)$ for some $\ctrv\in \mbX_h(\Sigma)$, and thus 
$\langle \Ts^{-1}(\ctrq),\ctrp\rangle = \langle \ctrv,\ctrp\rangle = 0$ by the very
definition \eqref{PolarSpace}. We conclude that $\mbX_h(\Sigma)^\circ$  and
$\Ts(\mbX_h(\Sigma))$ are $\Ts^{-1}$-orthogonal to each other. \hfill $\Box$

\quad\\
The spaces $\mbX_h(\Sigma)$ and $\mbX_h(\Sigma)^\circ$ will play an
important role in the sequel, and we need a convenient way to characterize them.
This is our motivation for introducing an oblique (i.e. non-self-adjoint) counterpart of
the exchange operator considered in \cite[Cor.5.1]{claeys2019new} and \cite[Lem.6.2]{claeys2020robust}. 

\begin{lem}\label{ObliqueIsometry}\quad\\
  Under Assumption \ref{CoercivityImpedance}, the operator $\Pi := (\mT+\mT^*)\mR(\mR^*\mT^*\mR)^{-1}\mR^*-\Id$
  with mapping property $\Pi:\mbVh(\Sigma)^*\to \mbVh(\Sigma)^*$ is an isometry in
  the norm \eqref{DefScalarProd} i.e. for all $\ctrp\in \mbVh(\Sigma)^*$ we have 
  \begin{equation*}
    \Vert \Pi(\ctrp)\Vert_{\Ts^{-1}} = \Vert \ctrp\Vert_{\Ts^{-1}}.
  \end{equation*}
\end{lem}
\noindent \textbf{Proof:}

Pick an arbitrary $\ctrp\in \mbVh(\Sigma)^*$ and set $\ctru := \mR(\mR^*\mT^*\mR)^{-1}\mR^*\ctrp$.
We have $\ctru \in\Range(\mR) =  \mbX_h(\Sigma)$ according to \eqref{Polarity} and
$\mR^*\mT^*(\bu) = \mR^*\bp$ hence, pairing with an arbitrary $w\in\Vh(\Sigma)$,
we obtain $\langle \mT^*(\bu),\mR(w)\rangle = \langle\bp,\mR(w)\rangle$.
Since $\mR:\Vh(\Sigma)\to \mbX_h(\Sigma)$ is surjective and $w$ is chosen arbitrarily in $\Vh(\Sigma)$,
we conclude that $\bu$ is characterized as the unique solution to the following variationnal problem
\begin{equation}
  \begin{aligned}
    & \ctru\in \mbX_h(\Sigma)\quad \text{and}\\
    & \langle \mT^*(\ctru), \ctrv\rangle = \langle \ctrp,\ctrv\rangle \quad\forall \ctrv\in \mbX_h(\Sigma).
  \end{aligned}
\end{equation}
In particular, taking $\ctrv = \overline{\ctru}$, we obtain
\begin{equation}
  \begin{aligned}
    2 \Re e\{\langle \ctrp,\overline{\ctru}\rangle\}
    & = 2\Re e\{\langle \mT^*(\ctru),\overline{\ctru}\rangle\}
    = \langle \mT^*(\ctru),\overline{\ctru}\rangle + \overline{\langle\overline{\ctru}, \mT^*(\ctru)\rangle}\\
    & = \langle \mT^*(\ctru),\overline{\ctru}\rangle + \langle\ctru, \overline{\mT^*(\ctru)}\rangle
    = 2\langle \Ts(\ctru),\overline{\ctru}\rangle\\
    & = 2\Vert \ctru\Vert_{\Ts}^2 = 2\Vert \Ts(\ctru)\Vert_{\Ts^{-1}}^2.
  \end{aligned}
\end{equation}
Next, observe that $\Pi(\ctrp) = 2\Ts(\ctru)-\ctrp$ by construction, which leads to
\begin{equation*}
  \begin{aligned}
    \Vert \Pi(\ctrp)\Vert_{\Ts^{-1}}^{2}
    & = \Vert 2\Ts(\ctru)-\ctrp \Vert_{\Ts^{-1}}^{2}\\
    & = \Vert \ctrp\Vert_{\Ts^{-1}}^{2} - 4\Re e\{\langle \ctrp,\overline{\ctru}\rangle\} +4\Vert\Ts(\ctru)\Vert_{\Ts^{-1}}^{2}
    = \Vert \ctrp\Vert_{\Ts^{-1}}^{2}.
  \end{aligned}
\end{equation*}
\hfill $\Box$

\quad\\
Interestingly, since $\Pi$ is a $\Ts^{-1}$-isometry and $\mbVh(\Sigma)^*$ is finite dimensional,
we deduce an expression of the inverse $\Pi^{-1} = (\mT+\mT^*)\mR(\mR^*\mT\mR)^{-1}\mR^*-\Id$.
This formula shall not be of much use in the present contribution though. In the special case 
of a self-adjoint impedance, this exchange operator becomes an orthogonal symmetry.

\begin{lem}\label{ProprieteEchangeOrthogonal}\quad\\
  Let Assumption \ref{CoercivityImpedance} hold and suppose further
  that $\mT = \mT^*$. Then  $\Pi = 2\mT \mR(\mR^*\mT\mR)^{-1}\mR^*-\Id$ defined
  in Lemma \ref{ObliqueIsometry} satisfies in addition the
  following properties
  \begin{itemize}
  \item[i)] $(\Id-\Pi)/2$ is the $\mT^{-1}$-orthogonal
    projection onto $\mbX_h(\Sigma)^\circ$,\\[-15pt]
  \item[ii)] $(\Id+\Pi^*)/2$ is the $\mT$-orthogonal
    projection onto $\mbX_h(\Sigma)$,\\[-15pt]
  \item[iii)]   $\Pi^2 = \Id$,\\[-15pt]
  \item[iv)] $\Pi^*\mR = \mR$.
  \end{itemize}
\end{lem}
\noindent \textbf{Proof:}

In the case $\mT = \mT^*$ we have $\mT = \Ts$. 
The identities \textit{iii)} and \textit{iv)} follow from direct calculus.
Then it is clear that $(\Id\pm\Pi)/2$ and  $(\Id\pm\Pi^*)/2$ are projectors.
The $\mT^{-1}$-orthogonality of $(\Id-\Pi)/2$ and the $\mT$-orthogonality of 
$(\Id+\Pi^*)/2$ are equivalent to the  $\mT^{-1}$-orthogonality of
$(\Id+\Pi)/2 =\mT \mR(\mR^*\mT\mR)^{-1}\mR^*$. The latter is proved by
taking $\ctrp,\ctrq\in\mbVh(\Sigma)^*$ arbitrary, and observing that
\begin{equation*}
  \begin{aligned}
    \langle \mT^{-1}(\Id+\Pi)\ctrp,\overline{\ctrq}\rangle 
    & = 2\langle (\mR^*\mT\mR)^{-1}\mR^*\ctrp, \mR^*\overline{\ctrq}\rangle \\
    & = 2\langle \ctrp, \mR(\mR^*\mT\mR)^{-1}\mR^*\overline{\ctrq}\rangle
    = \langle \mT^{-1}\ctrp, \overline{(\Id+\Pi)\ctrq}\rangle
  \end{aligned}
\end{equation*}
\hfill $\Box$

\quad\\
The previous lemma is consistent with the theory proposed in our previous contributions
that only considered hermitian positive definite (HPD) impedances, see Corollary 5.1 in \cite{claeys2019new}
and Lemma 6.2 in \cite{claeys2020robust}. Coming back to the general case of a priori non-self-adjoint
impedance, let us show how this exchange operator can serve for the effective characterization of
$\mbX_h(\Sigma)\times \mbX_h(\Sigma)^\circ$.

\begin{lem}\label{CaracTransCond}\quad\\
  Let Assumption \ref{CoercivityImpedance} hold and 
  define $\mathscr{X}_h(\Sigma):=\mbX_h(\Sigma)\times \mbX_h(\Sigma)^\circ$ and
  $\Pi := 2\Ts\mR(\mR^*\mT^*\mR)^{-1}\mR^*-\Id$ as in Lemma \ref{ObliqueIsometry}.
  Then for any pair $(\ctru,\ctrp)\in \mbVh(\Sigma)\times\mbVh(\Sigma)^*$
  we have the equivalence
  \begin{equation}\label{CaracTransmissionConditions}
    (\ctru,\ctrp)\in \mathscr{X}_h(\Sigma)
    \quad\iff\quad -\ctrp + i\mT(\ctru) = \Pi(\ctrp + i\mT^*(\ctru)).
  \end{equation}
\end{lem}
\noindent \textbf{Proof:}

Let us first investigate the action of $\Pi$ on certain well chosen input arguments.
If $\ctrp\in \mbX_h(\Sigma)^\circ$ then $\mR^*(\ctrp) = 0$ according to \eqref{Polarity}
hence $\Pi(\ctrp) = -\ctrp$. On the other hand, if $\ctru\in \mbX_h(\Sigma)$ then  $\ctru = \mR(v)$
for some $v\in \Vh(\Sigma)$, and this yields
\begin{equation}
  \begin{aligned}
    \Pi\cdot\mT^*(\ctru)
    & =  (\mT+\mT^*)\mR(\mR^*\mT^*\mR)^{-1}(\mR^*\mT^*\mR)v-\mT^*(\ctru)\\
    & =  (\mT+\mT^*)\mR v -\mT^*(\ctru)= \mT(\ctru).
  \end{aligned}
\end{equation}
The previous observations clearly imply that any pair  $(\ctru,\ctrp)\in
\mbX_h(\Sigma)\times \mbX_h(\Sigma)^\circ$ satisfies
$-\ctrp + i\mT(\ctru) = \Pi(\ctrp + i\mT^*(\ctru))$. Reciprocally pick an arbitrary
pair $(\ctru,\ctrp)\in \mbVh(\Sigma)\times \mbVh(\Sigma)^*$ satisfying
$-\ctrp + i\mT(\ctru) = \Pi(\ctrp + i\mT^*(\ctru))$. Plugging the expression of $\Pi$
into this identity implies
\begin{equation}\label{CaracIntermediaire}
  \begin{aligned}
    && -\ctrp + i\mT(\ctru) & =  (\mT+\mT^*)\mR(\mR^*\mT^*\mR)^{-1}\mR^*(\ctrp + i\mT^*(\ctru))- \ctrp - i\mT^*(\ctru)\\
    \iff&& i(\mT+\mT^*)\ctru & = (\mT+\mT^*)\mR(\mR^*\mT^*\mR)^{-1}\mR^*(\ctrp+i\mT^*(\ctru))\\
  \end{aligned}
\end{equation}
Next multiply the last identity above on the left by $\mR^*\mT^*(\mT+\mT^*)^{-1}$
which leads to $i\mR^*\mT^*(\ctru) =  i\mR^*\mT^*(\ctru) + \mR^*\ctrp$ and finally
$\mR^*\ctrp = 0$ which is equivalent to $\ctrp\in \mbX_h(\Sigma)^\circ$
according to \eqref{Polarity}. Coming back to the second identity of
\eqref{CaracIntermediaire}, and multiplying by $-i(\mT+\mT^*)^{-1}$ and taking
account of \eqref{Polarity}, we obtain $\ctru = \mR(\mR^*\mT^*\mR)^{-1}\mR^*\mT^*(\ctru)\in
\mrm{Range}(\mR) = \mbX_h(\Sigma)$. This finishes the proof. \hfill $\Box$

\begin{remark}\label{NonInvolutive}
The exchange operator $\Pi$ defined in Lemma \ref{CaracTransCond}
  is not a priori involutive i.e. $\Pi^2\neq \Id$ in general. Such a property is garanteed only
  when the impedance operator $\mT$ is HPD, see \textit{iii)} of Lemma
  \ref{ProprieteEchangeOrthogonal}.

To illustrate this, consider Example \ref{ExampleImpedancezTR} and let us examine the
  exchange operator in this case. We have $\mT^* = \overline{z}\mT_{\textsc{r}}$ and
  $\mT+\mT^* = (z+\overline{z}) \mT_{\textsc{r}}$. The operator $\mP_{\textsc{r}}:=
  \mT_{\textsc{r}}\mR(\mR^*\mT_{\textsc{r}}\mR)^{-1}\mR^{*}$ is then a $\Ts^{-1}-$orthogonal
  projector, and a direct calculation yields
  \begin{equation*}
    \begin{aligned}
      \Pi^{\textcolor{white}{2}}
      & = \frac{z+\overline{z}}{\overline{z}}\mP_{\textsc{r}}-\Id
      = \Big(\frac{z}{\vert z\vert}\Big)^{2}\mP_{\textsc{r}}-(\Id-\mP_{\textsc{r}})\\
      \Pi^2
      & =  \Big(\frac{z}{\vert z\vert}\Big)^4\mP_{\textsc{r}}+(\Id-\mP_{\textsc{r}})
    \end{aligned}
  \end{equation*}
  We see that $\Pi^2 = \Id\iff (z/\vert z\vert)^4 = 1 \iff z\in \RR\cup i\RR$.
  Taking account of Assumption \ref{CoercivityImpedance}, $\Pi$ is an involution
  only if $z\in (0,+\infty)$, in which case $\mT$ is HPD.
\end{remark}

\section{Reformulation of the scattering problem}
The equivalence \eqref{CaracTransmissionConditions} is a new way to impose transmission
conditions across interfaces between subdomains. We now make use of this new ingredient
to reformulate the wave propagation problem we are interested in. Let us introduce an
operator  $\Ah:\mbVh(\Omega)\to \mbVh(\Omega)^*$ and a source term  $\ctrf\in \mbVh(\Omega)^*$
associated to the domain decomposed problem and defined by
\begin{equation}\label{DefOperators}
  \begin{array}{l}
    \langle \Ah\ctru,\ctrv\rangle:= \sum_{j=1\dots \mJ}\int_{\Omega_{j}}\mu\nabla u_j\cdot\nabla v_j-\kappa^{2} u_jv_j \;d\bx\\[5pt]
    \langle \ctrf,\ctrv\rangle:= \sum_{j=1\dots \mJ} \int_{\Omega_{j}}f v_j d\bx + \int_{\partial \Omega_{j}\cap\partial \Omega }g v_j d\sigma
  \end{array}
\end{equation}
for any $\ctru = (u_1,\dots,u_\mJ), \ctrv = (v_1,\dots,v_\mJ)$ in $\mbVh(\Omega)$. These are nothing but
a finite element matrix and vector. A few remarks are in order concerning the operator $\Ah$. First
of all it admits a block diagonal form with respect to the subdomain decomposition. Taking account
of Assumption \ref{Hypo1}, the imaginary part of $\Ah$ is signed
\begin{equation}\label{Absorption}
  \Im m\{ \langle \Ah \ctru,\overline{\ctru}\rangle\}\leq 0\quad \forall \ctru\in \mbVh(\Omega).
\end{equation}
Besides, a comparison of \eqref{DefOperators} with \eqref{BilinearFormPb}-\eqref{ContinuousNorms} shows  that,
if $\ctru = (u\vert_{\Omega_1},\dots,u\vert_{\Omega_\mJ})$ and $\ctrv =  (v\vert_{\Omega_1},\dots,v\vert_{\Omega_\mJ})$
for some $u,v\in \Vh(\Omega)$, then  $\ctru,\ctrv\in \mbX_h(\Omega)$ according to \eqref{ContinuousNorms}
and, in this case, we have $\langle \Ah\ctru,\overline{\ctrv}\rangle = a(u,v)$ and
$\langle \ctrf,\overline{\ctrv}\rangle = \ell(v)$. We introduce the continuity modulus $\Vert a\Vert$
of $\Ah$, and also re-express the inf-sup constant $\alpha_h$ from Assumption \ref{Hypo2} as follows
\begin{equation}\label{MatrixInfSup}
  \begin{aligned}
    \alpha_h = & \mathop{\inf\phantom{p}}_{\ctru\in \mbX_h(\Omega)\setminus\{0\}}
    \sup_{\ctrv\in\mbX_h(\Omega)\setminus\{0\}}\frac{\vert \langle \Ah(\ctru),
      \ctrv\rangle\vert}{\Vert \ctru\Vert_{\mbH^1(\Omega)}\Vert \ctrv\Vert_{\mbH^1(\Omega)}}>0\\
    \Vert a \Vert := & \sup_{\ctru\in \mbVh(\Omega)\setminus\{0\}}
    \sup_{\ctrv\in\mbVh(\Omega)\setminus\{0\}}\frac{\vert \langle \Ah(\ctru),
      \ctrv\rangle\vert}{\Vert \ctru\Vert_{\mbH^1(\Omega)}\Vert \ctrv\Vert_{\mbH^1(\Omega)}}
  \end{aligned}
\end{equation}
The constant $\Vert a\Vert$ is simply the continuity modulus of the sesquilinear form
$a(\cdot,\cdot)$ over the space of piecewise $\mathbb{P}_k$-Lagrange functions.
We will now re-write Problem \eqref{DiscreteVF} in several equivalent forms more prone
to domain decomposition. First we use operator $\Ah$ so as to put it in matrix form.
\begin{lem}\label{LemmaReformulation1}\quad\\
  Assume that $u\in \Vh(\Omega)$ is solution to \eqref{DiscreteVF}. Then, setting
  $\ctru = (u\vert_{\Omega_1},\dots,u\vert_{\Omega_{\mJ}})\in \mbVh(\Omega)$, there exists
  $\ctrp\in \mbVh(\Sigma)^*$ such that
  \begin{equation}\label{MatrixEq1}
    \begin{aligned}
      & (\ctru,\ctrp)\in \mbX_h(\Omega)\times\mbX_h(\Sigma)^\circ \;\text{and}\\
      & \Ah\ctru-\Bh^*\ctrp = \ctrf.
    \end{aligned}
  \end{equation}
  Reciprocally if the pair $(\ctru,\ctrp)\in \mbVh(\Omega)\times \mbVh(\Sigma)^*$ solves \eqref{MatrixEq1},
  with $\ctru = (u_1,\dots,u_\mJ)$,  then the function defined by $u(\bx) = u_1(\bx)1_{\Omega_{1}}(\bx) + \dots + u_\mJ(\bx)1_{\Omega_{\mJ}}(\bx)$
  belongs to $\Vh(\Omega)$ and solves \eqref{DiscreteVF}. 
\end{lem}
\noindent \textbf{Proof:}

Assume that $u\in \Vh(\Omega)$ is solution to \eqref{DiscreteVF}. By definition we have
$\ctru\in \mbX_h(\Omega)$ where $\ctru = (u\vert_{\Omega_1},\dots,u\vert_{\Omega_{\mJ}})$
and \eqref{DiscreteVF} rewrites
\begin{equation}\label{Reformulation1}
  \begin{aligned}
    &\ctru \in \mbX_h(\Omega)\;\;\text{such that}\\
    &\langle \Ah(\ctru),\ctrv\rangle = \langle \ctrf,\ctrv\rangle\;\;
    \forall \ctrv\in  \mbX_h(\Omega)
  \end{aligned}
\end{equation}
Since $\bv\in \mbX_h(\Omega)$ whenever $\bv\in \mbVh(\Omega)$ satisfies $\Bh(\bv) = 0$, we can apply
Lemma \ref{LiftingBoundaryFunctional} to $\Ah(\ctru) - \ctrf\in \mbVh(\Omega)^*$ which yields the
existence of $\ctrp \in \mbVh(\Sigma)^*$ such that $\Ah\ctru - \ctrf = \Bh^*\ctrp$ which is the second line
of \eqref{MatrixEq1}. Then the second line of \eqref{Reformulation1} rewrites $0 = \langle \Bh^*\ctrp,\ctrv\rangle
= \langle \ctrp,\Bh(\ctrv)\rangle\;\forall \bv\in \mbX_h(\Omega)$. Since $\Bh$ maps $\mbX_h(\Omega)$ onto
$\mbX_h(\Sigma)$, we conclude that $\ctrp\in \mbX_h(\Sigma)^\circ$. As a consequence \eqref{MatrixEq1} holds. 

\quad\\
Now assume that $(\ctru,\ctrp)\in \mbVh(\Omega)\times \mbVh(\Sigma)^*$ solves \eqref{MatrixEq1},
where we denote $\ctru = (u_1,\dots,u_\mJ)$.  Besides, since $\ctrp\in\mbX_h(\Sigma)^\circ$, we have
$\langle \Bh^*\ctrp,\ctrv\rangle = \langle\ctrp, \Bh(\ctrv)\rangle=0$ for all $\ctrv\in \mbX_h(\Omega)$
which rewrites as \eqref{Reformulation1}. Since \eqref{Reformulation1} implies  \eqref{DiscreteVF}
with the function $u(\bx) = u_1(\bx)1_{\Omega_{1}}(\bx) + \dots + u_\mJ(\bx)1_{\Omega_{\mJ}}(\bx)$ belonging to $\Vh(\Omega)$,
this concludes the proof. \hfill $\Box$

\quad\\
The tuple of unkowns $\ctrp$ in Formulation \eqref{MatrixEq1} should be understood as Neumann fluxes
of the volume solution across boundaries of subdomains, and the condition $\ctrp\in\mbX_h(\Sigma)^\circ$
should be understood as the Neumann part of classical transmission conditions across interfaces. 

Note that $\ctru\in \mbX_h(\Omega)$ if and only if $\ctru\in \mbVh(\Omega)$ and $\Bh(\ctru)\in \mbX_h(\Sigma)$. 
The property  $\Bh(\ctru)\in \mbX_h(\Sigma)$ should be interpreted as the Dirichlet part of classical transmission
conditions across interfaces. Thanks to the previous remarks, we can transform further Formulation \eqref{MatrixEq1}
by taking account of the characterization of $\mbX_h(\Sigma)\times \mbX_h(\Sigma)^\circ$ stemming from Lemma \ref{CaracTransCond}.
This directly yields the following reformulation.

\begin{lem}\label{Reformulation2}\quad\\
  Under Assumption \ref{CoercivityImpedance}, the pair $(\ctru,\ctrp)$ solves \eqref{MatrixEq1}
  if and only if it satisfies
  \begin{equation}\label{MatrixEq2}
    \begin{aligned}
      & (\ctru,\ctrp)\in \mbVh(\Omega)\times \mbVh(\Sigma)^*\;\text{and}\\
      & \Ah\ctru-\Bh^*\ctrp = \ctrf,\\
      & \ctrp - i\mT\Bh(\ctru) = -\Pi(\ctrp + i\mT^*\Bh(\ctru)).
    \end{aligned}
  \end{equation}
\end{lem}

\noindent 
The formulation above only involves the multi-trace space and its dual where contributions from subdomains
are decorrelated from each other. Transmission conditions across interfaces come only into play through the
exchange operator $\Pi$. Although we do not write it as an iterative algorithm, Formulation \eqref{MatrixEq2}
above is entirely analogous to e.g. \cite[Eq.(8)]{despres:hal-03230250}, \cite[Eq.(24)]{claeys2020robust},
\cite[Eq.(3)-(4)]{Boubendir2012} or \cite[Eq.(59)]{MR1764190}.

\section{Scattering operator}

The next step of our analysis consists in eliminating the volume unknowns in \eqref{MatrixEq2}.
However the map $\Ah:\mbVh(\Omega)\to \mbVh(\Omega)^*$ may be not invertible, so we avoid
using $\Ah^{-1}$ and re-arrange local subproblems. The next lemma establishes that
local subproblems become invertible if we add absorbing impedance conditions.
\begin{lem}\quad\\
  Under Assumptions \ref{Hypo1}, \ref{Hypo2} and \ref{CoercivityImpedance},
  the operator $\Ah-i\Bh^*\mT\Bh:\mbVh(\Omega)\to \mbVh(\Omega)^*$ is systematically invertible
  which rewrites in inf-sup condition form as
  \begin{equation}
    \beta_h:= \mathop{\inf\phantom{p}}_{\ctru\in \mbVh(\Omega)\setminus\{0\}}
    \sup_{\ctrv\in \mbVh(\Omega)\setminus\{0\}} \frac{\vert \langle (\Ah-i\Bh^*\mT\Bh)\ctru,\ctrv
      \rangle\vert}{\Vert \ctru\Vert_{\mbH^1(\Omega)}\Vert \ctrv\Vert_{\mbH^1(\Omega)}}>0.
  \end{equation}
\end{lem}
\noindent \textbf{Proof:}

Assume that $\beta_h = 0$ so that there exists $\ctru\in \mbVh(\Omega)$ such that 
$(\Ah-i\Bh^*\mT\Bh)\ctru = 0$. From this and \eqref{Absorption}, we conclude
that $0 =-\Im m \{\langle \Ah(\ctru),\overline{\ctru}\rangle\} +
\Im m \{i\langle \mT\Bh(\ctru),\Bh(\overline{\ctru})\rangle\}\geq
\Re e\{\langle \mT\Bh(\ctru),\Bh(\overline{\ctru})\rangle\}$
which implies $\Bh(\ctru) = 0$ according to Assumption \ref{CoercivityImpedance}.
Since we have $\mrm{Ker}(\Bh)\subset \mbX_h(\Omega)$ according to \eqref{Polarity},
this shows that $\ctru\in \mbX_h(\Omega)$ and $\Ah(\ctru) = 0$,
which implies $\ctru = 0$ according to \eqref{MatrixInfSup}. This establishes that
$\mrm{Ker}(\Ah-i\Bh^*\mT\Bh) = \{0\}$ and finishes the proof. \hfill $\Box$

\quad\\
Now we can introduce a so called scattering operator. 
The following result is inspired by \cite[Lemma 6]{despres:hal-02612368}. 

\begin{lem}\label{ScatteringOperator}\quad\\
  Under Assumptions \ref{Hypo1}, \ref{Hypo2} and \ref{CoercivityImpedance},
  the operator $\Sh := \Id + 2i\Ts\Bh (\Ah-i\Bh^*\mT\Bh)^{-1}\Bh^*$
  is a $\Ts^{-1}$-contraction and, for all $\ctrp\in \mbVh(\Sigma)^*$,
  satisfies the identity
  \begin{equation}\label{EnergyDecay}
    \begin{aligned}
      & \Vert \Sh(\ctrp)\Vert_{\Ts^{-1}}^2 + 4\vert \Im m\{\langle \Ah\ctru,\overline{\ctru}\rangle\}\vert
      = \Vert \ctrp\Vert_{\Ts^{-1}}^2\\
      & \text{where}\quad \ctru = (\Ah-i\Bh^*\mT\Bh)^{-1}\Bh^*\ctrp. 
    \end{aligned}
  \end{equation}
\end{lem}
\noindent \textbf{Proof:}

We note that $\Sh(\ctrp) = \ctrp +2i\Ts\Bh(\ctru)$, and simply expand the left hand side of \eqref{EnergyDecay},
taking account of the sign property provided by \eqref{Absorption}. This yields the following calculus
\begin{equation}
  \begin{aligned}
    \Vert \Sh(\ctrp)\Vert_{\Ts^{-1}}^2
    & = \Vert \ctrp +2i\Ts\Bh(\ctru) \Vert_{\Ts^{-1}}^2\\
    & = \Vert \ctrp\Vert_{\Ts^{-1}}^2 -4 \Re e\{ i\langle \ctrp, \Bh(\overline{\ctru})\rangle\} + 4\Vert \Ts\Bh(\ctru) \Vert_{\Ts^{-1}}^2\\
    & = \Vert \ctrp\Vert_{\Ts^{-1}}^2 -4 \Re e\{ i\langle (\Ah-i\Bh^*\mT\Bh)\ctru, \overline{\ctru}\rangle\}
    + 4\langle \Ts\Bh(\ctru),\Bh(\overline{\ctru})\rangle\\
    & = \Vert \ctrp\Vert_{\Ts^{-1}}^2 +4\Im m\{\langle \Ah(\ctru),\overline{\ctru}\rangle\}
    = \Vert \ctrp\Vert_{\Ts^{-1}}^2 -4\vert \Im m\{\langle \Ah(\ctru),\overline{\ctru}\rangle\}\vert
  \end{aligned}
\end{equation}
\hfill $\Box$

\quad\\
The scattering operator $\mS$ is subdomain-wise block diagonal under the additional assumption that  $\mT$
is subdomain-wise block-diagonal, which is a reasonnable and easy property to fulfill in practice.

Identity \eqref{EnergyDecay} should be interpreted as energy conservation. It is similar to
\cite[Lemma 4.1]{MR1291197}, \cite[Lemma 1]{MR1764190} or \cite[Lemma 11]{despres:hal-03230250}. 
The scattering operator $\Sh$ takes a Neuman multi-trace as input, solves in each subdomain the
associated ingoing impedance problem, and returns the outgoing impedance trace
as output. This interpretation of the scattering operator is made explicit in the next result.

\begin{prop}\label{InterpretationScattering}\quad\\
  Define $\mathscr{C}_h(\Sigma):=\{(\ctrv,\ctrp)\in \mbVh(\Sigma)\times\mbVh(\Sigma)^*:
  \;\exists \ctru\in \mbVh(\Omega),\;\Ah\ctru = \Bh^*\ctrp,\;\Bh\ctru = \ctrv\}$
  which will be called the space of discrete Cauchy data. Let
  Assumptions \ref{Hypo1}, \ref{Hypo2} and \ref{CoercivityImpedance} hold. 
  Then for any pair $(\ctrv,\ctrp)\in \mbVh(\Sigma)\times\mbVh(\Sigma)^*$,
  we have the equivalence
  \begin{equation}\label{CauchyDataCaract}
    (\ctrv,\ctrp)\in\mathscr{C}_h(\Sigma)\quad \iff \quad \ctrp+i\mT^*\ctrv = \Sh(\ctrp-i\mT\ctrv).
  \end{equation}
\end{prop}
\noindent \textbf{Proof:}

Take an arbitrary pair $(\ctrv,\ctrp)\in\mathscr{C}_h(\Sigma)$ and set $\ctrq = \ctrp-i\mT\ctrv$.
By definition of the Cauchy data space, there exists $\ctru\in \mbVh(\Omega)$ such that
$\Ah\ctru = \Bh^*\ctrp$ and $\Bh\ctru = \ctrv$ which implies  in particular that
$(\Ah - i\Bh^*\mT\Bh)\ctru= \Bh^*\ctrq$ hence $\ctru = (\Ah - i\Bh^*\mT\Bh)^{-1}\Bh^*\ctrq$.
Next the definition of the scattering operator given in Lemma \eqref{ScatteringOperator} yields
$\Sh(\ctrq) = \ctrp-i\mT\ctrv + 2i\Ts\Bh(\ctru) = \ctrp-i\mT\ctrv + i\mT\ctrv + i\mT^*\ctrv = \ctrp
+ i\mT^*\ctrv$ which rewrites as \eqref{CauchyDataCaract}.

Reciprocally take any $(\ctrv,\ctrp)\in\mbVh(\Sigma)\times\mbVh(\Sigma)^*$ satisfying
$\ctrp+i\mT^*\ctrv = \Sh(\ctrp-i\mT\ctrv)$. Define $\ctru\in \mbVh(\Omega)$ by $\ctru = (\Ah - i\Bh^*\mT\Bh)^{-1}\Bh^*(\ctrp -i\mT\ctrv)$
so that, using the definition of $\Sh$ and the invertibility of $\Ts := (\mT + \mT^*)/2$, we have
$\ctrp+i\mT^*\ctrv = \Sh(\ctrp-i\mT\ctrv) = \ctrp-i\mT\ctrv + 2i\Ts\Bh\ctru \iff 2i\Ts\ctrv = 2i\Ts\Bh\ctru \iff \ctrv = \Bh\ctru$. 
Next coming back to the definition of $\ctru$ we conclude that $\Ah\ctru - i\Bh^*\mT\ctrv = \Bh^*(\ctrp -i\mT\ctrv)$,
and finally $\Ah\ctru = \Bh^*\ctrp$. This proves that the pair $(\ctrv,\ctrp)$ belongs to $\mathscr{C}_h(\Sigma)$. 
\hfill $\Box$

\quad\\
We draw the attention of the reader on the similarities between \eqref{CauchyDataCaract} and \eqref{CaracTransmissionConditions}.
Both characterizations are expressed in terms of ingoing traces (i.e. traces of the form $\ctrp-i\mT\ctrv$) and outgoing traces
(i.e. traces of the form  $\ctrp+i\mT^*\ctrv$). Let us also point an interesting property satisfied by the elements of $\mathscr{C}_h(\Sigma)$.

\begin{lem}\label{EnergyConservationSelfAdjoint}\quad\\  
  Under Assumptions \ref{Hypo1}, \ref{Hypo2} and \ref{CoercivityImpedance},
  if $(\ctrv,\ctrp)\in \mathscr{C}_h(\Sigma)$ then $\Im m\{\langle \ctrp,\overline{\ctrv}\rangle\}\leq 0$
  and we have the energy identities
  \begin{equation*}
    \begin{aligned}
      \Vert \ctrv\Vert_{\Ts}^2 + \Vert \ctrp\Vert_{\Ts^{-1}}^2
      & = \Vert \ctrp+i\Ts\ctrv\Vert_{\Ts^{-1}}^2 + 2\vert \Im m\{\langle \ctrp,\overline{\ctrv}\rangle\}\vert\\
      & = \Vert \ctrp-i\Ts\ctrv\Vert_{\Ts^{-1}}^2 - 2\vert \Im m\{\langle \ctrp,\overline{\ctrv}\rangle\}\vert.
    \end{aligned}
  \end{equation*}
\end{lem}
\noindent \textbf{Proof:}

By the very definition of $\mathscr{C}_h(\Sigma)$ given in Proposition \ref{InterpretationScattering},
there exists $\ctru\in \mbVh(\Omega)$ such that $\Ah(\ctru) = \Bh^*(\ctrp)$ and $\Bh(\ctru) = \ctrv$,
hence $\langle \ctrp,\overline{\ctrv}\rangle = \langle \ctrp,\Bh(\overline{\ctru})\rangle =
\langle \Bh^*(\ctrp),\overline{\ctru}\rangle = \langle \Ah(\ctru),\overline{\ctru}\rangle$.
With \eqref{Absorption} we obtain $\Im m\{\langle \ctrp,\overline{\ctrv}\rangle\} =
\Im m\{\langle \Ah(\ctru),\overline{\ctru}\rangle\}\leq 0$ which is the first desired
result. The energy identities directly follow from
$\Vert \ctrp\pm i\Ts\ctrv\Vert_{\Ts^{-1}}^2 =   \Vert \ctrp\Vert_{\Ts^{-1}}^2 \pm 2
\Im m\{\langle \ctrp,\overline{\ctrv}\rangle\} + \Vert \ctrv\Vert_{\Ts}^2$.
\hfill $\Box$

\section{Skeleton formulation}

We will now use the scattering operator introduced in the previous section
to rewrite equivalently our wave propagation problem as an equation posed on
the skeleton of the subdomain partition. In spite of a possible minor difference due to
a sign convention, the formulation derived below is
similar to \cite[Eq.(27)]{claeys2020robust}, \cite[Eq.(7.2)]{claeys2019new},
\cite[Eq.(45)\&(51)]{MR1764190}, \cite[\S3.3]{MR1291197},\cite[Chap.6]{lecouvez:tel-01444540},
\cite{despres:hal-03230250}.

\begin{lem}\label{EquivalenceSkeletonFormulation}\quad\\
  Let Assumptions \ref{Hypo1}, \ref{Hypo2} and \ref{CoercivityImpedance} hold.
  Set $\ctrg:=-2i\Pi\Ts\Bh(\Ah-i\Bh^*\mT\Bh)^{-1}\ctrf\in \mbVh(\Sigma)^*$. If the pair
  $(\ctru,\ctrp)\in \mbVh(\Omega)\times \mbVh(\Sigma)^*$ solves \eqref{MatrixEq2}
  then the tuple of traces $\ctrq = \ctrp - i\mT\Bh\ctru$ solves
  \begin{equation}\label{SkeletonEquation}
    \begin{aligned}
      & \ctrq\in \mbVh(\Sigma)^*\;\; \text{and}\\
      & (\Id +\Pi\Sh)\ctrq = \ctrg.
    \end{aligned}
  \end{equation}
  Reciprocally, if $\ctrq$ satisfies \eqref{SkeletonEquation} then the pair $(\ctru,\ctrp)$ defined by
  $\ctru = (\Ah-i\Bh^*\mT\Bh)^{-1}(\Bh^*\ctrq + \ctrf)$ and $\ctrp = \ctrq + i\mT\Bh\ctru$ solves \eqref{MatrixEq2}. 
\end{lem}
\noindent \textbf{Proof:}

Assume that $(\ctru,\ctrp)$ is solution to \eqref{MatrixEq2} and consider $\ctrq = \ctrp - i\mT\Bh\ctru$.
We have $\ctrp = \ctrq + i\mT\Bh\ctru$ and $\ctrp+i\mT^*\Bh\ctru = \ctrq + 2i\Ts\Bh\ctru$. Then we can replace
$\ctrp$ by $\ctrq + i\mT\Bh\ctru$ in \eqref{MatrixEq2} which leads to the equations 
\begin{equation*}
  \begin{aligned}
    & (\Ah - i\Bh^*\mT\Bh)\ctru = \Bh^*\ctrq + \ctrf,\\
    & \ctrq = -\Pi(\ctrq + 2i\Ts\Bh\ctru).
  \end{aligned}
\end{equation*}
Let us now decompose $\ctru = \tilde{\ctru} + \ctru_f$ where $\ctru_f = (\Ah - i\Bh^*\mT\Bh)^{-1}\ctrf$
and $\tilde{\ctru} = \ctru - \ctru_f$. This leads to  $(\Ah - i\Bh^*\mT\Bh)\tilde{\ctru} = \Bh^*\ctrq$ and
$\ctrq = -\Pi(\ctrq + 2i\Ts\Bh\tilde{\ctru})+\ctrg$. There only remains to eliminate $\tilde{\ctru}$ in these
equations, taking account of the definition of $\Sh$ from Lemma \ref{ScatteringOperator}, which yields 
$(\Id +\Pi\Sh)\ctrq = \ctrg$.

\quad\\
Reciprocally assume that $\ctrq\in \mbVh(\Sigma)^*$ solves \eqref{SkeletonEquation}, and set 
$\ctru = (\Ah-i\Bh^*\mT\Bh)^{-1}(\Bh^*\ctrq + \ctrf)$ and  $\ctrp = \ctrq + i\mT\Bh\ctru$.
As a consequence, by construction, the first equation of \eqref{MatrixEq2} is satisfied,
namely $\Ah\ctru = \Bh^*\ctrp + \ctrf$. There only remains to verify that the second equation of
\eqref{MatrixEq2} is satisfied by $(\ctru,\ctrp)$ as well.
Set $\ctru_f = (\Ah-i\Bh^*\mT\Bh)^{-1}\ctrf$ so that $\ctru-\ctru_f = (\Ah-i\Bh^*\mT\Bh)^{-1}\Bh^*\ctrq$
hence  $\Sh\ctrq = \ctrq + 2i\Ts\Bh(\ctru-\ctru_f)$. Plugging this in \eqref{SkeletonEquation} yields
\begin{equation*}
  \begin{aligned}
    \ctrq
    & = -\Pi\Sh(\ctrq) + \ctrg = -\Pi(\ctrq + 2i\Ts\Bh(\ctru-\ctru_f))+\ctrg\\
    & = -\Pi(\ctrq + 2i\Ts\Bh\ctru)
  \end{aligned}
\end{equation*}
Since $\ctrq = \ctrp - i\mT\Bh\ctru$ and $\ctrq + 2i\Ts\Bh\ctru = \ctrp+i\mT^*\Bh\ctru$,
the above equation rewrites $-\ctrp + i\mT\Bh\ctru = \Pi(\ctrp+i\mT^*\Bh\ctru)$ which concludes
the proof. \hfill $\Box$

\quad\\
Previously we exhibited a chain of equivalent formulations that relates
the initial discrete variational problem \eqref{DiscreteVF} to the skeleton
equation \eqref{SkeletonEquation}. Well posedness of \eqref{DiscreteVF} thus
immediately implies well posedness of\eqref{SkeletonEquation} but there is
more: the skeleton formulation provides a strongly coercive formulation of
our Helmholtz problem. The following result is completely similar to
Corollary 8.4 of \cite{claeys2020robust}.

\begin{cor}\quad\\
  Under Assumptions \ref{Hypo1}, \ref{Hypo2} and \ref{CoercivityImpedance},
  the operator $\Id+\Pi\Sh:\mbVh(\Sigma)^*\to \mbVh(\Sigma)^*$ is an isomorphism
  that is $\Ts^{-1}$-coercive. More precisely, for all $\ctrq\in \mbVh(\Sigma)^*$
  we have
  \begin{equation}\label{CoercivityProperty}
    \begin{aligned}
      & \Re e\{\langle (\Id +\Pi\Sh)\ctrq,\Ts^{-1}(\overline{\ctrq})\rangle\}\geq
      \frac{\gamma_h^2}{2}\Vert \ctrq\Vert_{\Ts^{-1}}^2\\
      & \text{where}\quad \gamma_h := \inf_{\ctrq\in \mbVh(\Sigma)^*\setminus\{0\}}
      \Vert (\Id +\Pi\Sh)\ctrq\Vert_{\Ts^{-1}}/\Vert \ctrq\Vert_{\Ts^{-1}}
    \end{aligned}
  \end{equation}
\end{cor}
\noindent \textbf{Proof:}

We first prove that $\mrm{Ker}(\Id +\Pi\Sh)=\{0\}$ which will show that $\Id +\Pi\Sh$
is an isomorphism (because $\mrm{dim}(\mbVh(\Sigma)^*)<+\infty$) and $\gamma_h>0$.
Assume that $(\Id +\Pi\Sh)\ctrq = 0$ for some $\ctrq\in \mbVh(\Sigma)^*$.
Set $\ctrg = 0$, $\ctrf = 0$, $\ctru = (\Ah-i\Bh^*\mT\Bh)^{-1}\Bh^*\ctrq$ and $\ctrp = \ctrq+i\mT\Bh\ctru$.
Applying Lemma \ref{EquivalenceSkeletonFormulation} we see that the pair $(\ctru,\ctrp)$  solves
\eqref{MatrixEq2} with $\ctrf = 0$. Next applying the equivalences given by Lemma \ref{Reformulation2}
and \eqref{LemmaReformulation1}, and using Assumption \ref{Hypo2} that yields existence and uniqueness of
the solution to these boundary value problems, we conclude that $(\ctru,\ctrp) = (0,0)$, hence $\ctrq = 0$.

\quad\\
Next, combining Lemma \ref{ObliqueIsometry} and \ref{ScatteringOperator},
we see that $\Pi\Sh$ is a contraction with respect to the norm \eqref{DefScalarProd}
induced by $\Ts^{-1}$. As a consequence, we obtain 
\begin{equation*}
  \begin{aligned}
    & \Vert  \ctrq\Vert_{\Ts^{-1}}^2 \geq \Vert \Pi\Sh\ctrq\Vert_{\Ts^{-1}}^{2}
    =  \Vert \ctrq - (\Id+\Pi\Sh)\ctrq\Vert_{\Ts^{-1}}^{2}\\
    & \phantom{\Vert  \ctrq\Vert_{\Ts^{-1}}^2 \geq  \Vert \Pi\Sh\ctrq\Vert_{\Ts^{-1}}^{2} }
    = \Vert  \ctrq\Vert_{\Ts^{-1}}^2-2\Re e\{\langle (\Id+\Pi\Sh)\ctrq,\Ts^{-1}(\overline{\ctrq})\rangle\}
    + \Vert (\Id+\Pi\Sh)\ctrq\Vert_{\Ts^{-1}}^{2}\\
    & \Re e\{\langle (\Id+\Pi\Sh)\ctrq,\Ts^{-1}(\overline{\ctrq})\rangle\}/\Vert  \ctrq\Vert_{\Ts^{-1}}^2
    \geq\frac{1}{2} \Vert (\Id+\Pi\Sh)\ctrq\Vert_{\Ts^{-1}}^{2}/\Vert  \ctrq\Vert_{\Ts^{-1}}^2\geq
    \gamma_h^2/2. 
  \end{aligned}
\end{equation*}
\hfill $\Box$

\quad\\
According to Lemma \ref{ObliqueIsometry} and \ref{ScatteringOperator}, we have 
$\Vert \Pi\mS(\ctrp)\Vert_{\Ts^{-1}}\leq \Vert \ctrp\Vert_{\Ts^{-1}}\;\forall\ctrp\in \mbVh(\Sigma)^*$.
Taking account of \eqref{CoercivityProperty} in addition, we conclude that the field of values
in the $\Ts^{-1}$-scalar product $\{ \langle (\Id +\Pi\mS)\bp,\Ts^{-1}(\overline{\bp})\rangle:
\Vert\bp\Vert_{\Ts^{-1}}=1, \bp\in \mbVh(\Sigma)^*\}$ is contained in $\{\lambda\in \CC: \vert \lambda-1\vert\leq 1,
\Re e\{\lambda\}\geq \gamma_h^2/2\}$. Combined with e.g. Elman estimate \cite{zbMATH03831185},
this readily yields an upper bound on the rate of convergence of GMRes \cite[Prop. 10.35]{zbMATH06587724}.
A similar remark holds for e.g. Richardson's algorithm see e.g. \cite[\S 3.5]{zbMATH06587724}.

\section{Bounds on the trace operator}\label{SpectralBoundsTraceOperator}

Before conducting a more detailed convergence analysis for the skeleton equation \eqref{SkeletonEquation},
we need to derive a few estimates related to the trace operator. The most natural constants
related to this operator are its continuity modulus $t_h^+$ and its inf-sup constant $t_h^-$ defined by
\begin{equation}\label{BoundsTraceOperator}
  \begin{aligned}
    &
    t_h^- := \mathop{\inf\phantom{p}}_{\ctrp\in \mbVh(\Sigma)^*\setminus\{0\}}
    \sup_{\ctrv\in \mbVh(\Omega)\setminus\{0\}}\frac{\vert \langle \Bh(\ctrv),\ctrp\rangle\vert}{\Vert \ctrv\Vert_{\mbH^1(\Omega)}\Vert \ctrp\Vert_{\Ts^{-1}}},\\
    &
    t_h^+ := \sup_{\ctrp\in \mbVh(\Sigma)^*\setminus\{0\}}
    \sup_{\ctrv\in \mbVh(\Omega)\setminus\{0\}}\frac{\vert \langle \Bh(\ctrv),\ctrp\rangle\vert}{\Vert \ctrv\Vert_{\mbH^1(\Omega)}\Vert \ctrp\Vert_{\Ts^{-1}}}.
  \end{aligned}
\end{equation}
We have in particular $\Vert \Bh(\ctrv)\Vert_{\Ts}\leq t_h^+\Vert \ctrv\Vert_{\mbH^1(\Omega)}$.
We provide an important alternative interpretation of $t_h^\pm$.  Let us introduce a map
$\Bh^\dagger:\mbVh(\Sigma)\to \mbVh(\Omega)$ defined as Moore-Penrose pseudo-inverse (see e.g. \cite[\S5.5.4]{zbMATH06159604} or
\cite[\S 2.6]{zbMATH05189473}) of the trace operator with respect to the volume norm
\begin{equation}\label{HarmonicLifting}
  \begin{aligned}
    & \Bh\Bh^\dagger = \Id\quad \text{and}\\
    & \Vert \Bh^\dagger(\ctrv)\Vert_{\mbH^1(\Omega)} = \inf\{\Vert \bu\Vert_{\mbH^1(\Omega)}: \bu\in\mbVh(\Omega),\,\Bh(\bu) = \ctrv\}.
  \end{aligned}
\end{equation}
By construction, for any $\ctrp\in \mbVh(\Sigma)$ and any $\ctru\in\mbVh(\Omega)$ such that
$\Bh(\ctru) = \ctrp$ we have $\Vert \Bh^\dagger(\ctrp)\Vert_{\mbH^1(\Omega)}\leq \Vert \ctru \Vert_{\mbH^1(\Omega)}$.
The map $\Bh^\dagger$ is a classical object of domain decomposition literature that is sometimes referred
to as discrete harmonic lifting, see \cite[\S 4.4]{MR2104179} and \cite[Def. 1.55]{MR3013465}. Because 
$\Bh$ is subdomain-wise block-diagonal, $\Bh^\dagger$ is itself block-diagonal.
The operator $\Bh^\dagger\Bh:\mbH^1(\Omega)\to\mbH^1(\Omega)$ is a projector that is orthogonal
with respect to the scalar product induced by \eqref{VolumeH1Norm}.

\begin{lem}\quad\\
Define $\Lambda:\mbVh(\Sigma)\to \mbVh(\Sigma)^*$ as the unique symetric positive definite operator
satisfying $\langle \Lambda(\ctrv),\overline{\ctrv}\rangle:= \Vert \Bh^\dagger(\ctrv)\Vert_{\mbH^1(\Omega)}^2$
and set $\Vert \ctrv\Vert_\Lambda^2:=\langle \Lambda(\ctrv),\overline{\ctrv}\rangle$. Then we have the identities 
\begin{equation}\label{BoundsTraceOperator2}
    t_h^- = \mathop{\inf\phantom{p}}_{\ctrv\in \mbVh(\Sigma)\setminus\{0\}}
    \frac{\Vert\ctrv\Vert_{\Ts}}{\Vert \ctrv\Vert_{\Lambda_{\phantom{s}}}},\quad\quad 
    t_h^+ = \sup_{\ctrv\in \mbVh(\Sigma)\setminus\{0\}}
        \frac{\Vert\ctrv\Vert_{\Ts}}{\Vert \ctrv\Vert_{\Lambda_{\phantom{s}}}}.
\end{equation}
\end{lem}
\noindent \textbf{Proof:}

We only prove the identity related to $t_h^-$ because the proof of the other identity follows a very
similar path. First for any $\ctrv\in \mbVh(\Omega)$, setting  $\tilde{\ctrv} = \Bh^\dagger\Bh(\ctrv)$,
we have $\Bh(\ctrv) = \Bh(\tilde{\ctrv})$ and $\Vert \tilde{\ctrv}\Vert_{\mbH^1(\Omega)}\leq \Vert \ctrv\Vert_{\mbH^1(\Omega)}$.
Since the trace operator $\Bh:\mbVh(\Omega)\to \mbVh(\Sigma)$ is surjective, and $\Bh\Bh^\dagger\Bh = \Bh$ we conclude 
\begin{equation*}
  \begin{aligned}
    \sup_{\ctrv\in \mbVh(\Omega)\setminus\{0\}}\frac{\vert \langle \Bh(\ctrv),\ctrp\rangle\vert}{\Vert \ctrv\Vert_{\mbH^1(\Omega)}\Vert \ctrp\Vert_{\Ts^{-1}}}
    & =\sup_{\ctrv\in \mbVh(\Omega)\setminus\{0\}}\frac{\vert \langle \Bh(\ctrv),\ctrp\rangle\vert}{\Vert
      \Bh^\dagger\Bh(\ctrv)\Vert_{\mbH^1(\Omega)}\Vert \ctrp\Vert_{\Ts^{-1}}}\\
    & = \sup_{\ctru\in \mbVh(\Sigma)\setminus\{0\}}\frac{\vert \langle \ctru,\ctrp\rangle\vert}{\Vert\ctru\Vert_{\Lambda}\Vert \ctrp\Vert_{\Ts^{-1}}}
    = \frac{\Vert \ctrp\Vert_{\Lambda^{-1}}}{\Vert \ctrp\Vert_{\Ts^{-1}}}\\
  \end{aligned}
\end{equation*}
for any $\ctrp\in \mbVh(\Sigma)^*$ with $\Vert \ctrp\Vert_{\Lambda^{-1}}^2:=\langle \Lambda^{-1}\ctrp,\overline{\ctrp}\rangle$.
The above calculus leads to the Rayleigh quotient expression $t_h^- = \inf_{\ctrp\in \mbVh(\Sigma)^*\setminus\{0\}}\Vert \ctrp\Vert_{\Lambda^{-1}}/
\Vert \ctrp\Vert_{\Ts^{-1}}$, which shows that $(t_h^-)^2$ can be characterized as an extremum of a generalized
eigenvalue problem
\begin{equation*}
  \begin{aligned}
    (t_h^-)^2
    & = \min\{\lambda\in \RR: \mrm{Ker}(\Lambda^{-1}-\lambda \Ts^{-1})\neq \{0\}\}\\
    & = \min\{\lambda\in \RR: \mrm{Ker}(\Ts-\lambda\Lambda)\neq \{0\}\} =
    \inf_{\ctrv\in \mbVh(\Sigma)^\setminus\{0\}}\Vert \ctrv\Vert_{\Ts}^2/
    \Vert \ctrv\Vert_{\Lambda}^2
  \end{aligned}
\end{equation*}
\hfill $\Box$

\quad\\
From the previous result, we see that $\Vert \Bh^\dagger(\ctrv)\Vert_{\mbH^1(\Omega)} = \Vert \ctrv\Vert_{\Lambda}
\leq (1/t_h^-)\Vert \ctrv\Vert_{\Ts}$ for all $\ctrv\in\mbVh(\Sigma)$. The constants $t_h^\pm$ thus appear as key
constants in the interplay between $\Bh$ and $\Ts$.

\quad\\
Since $\Bh$ and $\Bh^\dagger$ are subdomain-wise block-diagonal, the operator $\Lambda:\mbVh(\Sigma)\to  \mbVh(\Sigma)^*$
is itself subdomain-wise block-diagonal $\Lambda = \mrm{diag}_{j=1\dots \mJ}(\Lambda_j)$, i.e. there are operators
$\Lambda_j:\Vh(\Gamma_j)\to \Vh(\Gamma_j)^*, j=1\dots \mJ$ such that 
\begin{equation*}
  \langle \Lambda(\ctru),\ctrv\rangle =   \langle \Lambda_1(u_1),v_1\rangle+\dots+\langle \Lambda_\mJ(u_\mJ),v_\mJ\rangle
\end{equation*}
for $\ctru,\ctrv\in \mbVh(\Sigma)$ with $\ctru = (u_1,\dots,u_\mJ)$ and $\ctrv = (v_1,\dots,v_\mJ)$. 
The operator $\Lambda$ can actually serve as an effective choice of self-adjoint impedance,
in which case $t_h^\pm = 1$. In any case, taking $\Lambda$ as a reference scalar product over $\mbVh(\Sigma)$,
the constants $t_h^\pm$ are  then the extremal eigenvalues of the symetrized impedance operator $\Ts$.
In addition, if $\mT$ is also subdomain-wise block-diagonal $\mT = \mrm{diag}_{j=1\dots \mJ}(\mT_{j})$,
then we have
\begin{equation*}
  \begin{aligned}
    & (t_h^-)^2 = \min_{j=1\dots \mJ}\inf_{v\in \Vh(\Gamma_j)\setminus\{0\}}
    \Re e\{\langle \mrm{T}_{j}(v),\overline{v}\rangle/\Vert v\Vert_{\Lambda_j}^2\},\\
    & (t_h^+)^2 = \max_{j=1\dots \mJ}\sup_{v\in \Vh(\Gamma_j)\setminus\{0\}}
        \Re e\{\langle \mrm{T}_{j}(v),\overline{v}\rangle/\Vert v\Vert_{\Lambda_j}^2\}.
  \end{aligned}
\end{equation*}
that is $t_h^\pm$ can be interpreted as bounds on the (real part of the)
field of values of the $\mT_j$'s with respect to local reference norms 
induced by the $\Lambda_j$'s. To summarize and since only the case of subdomain-wise
block-diagonal impedance can be regarded as computationally reasonnable a situation,
in practice, the constants $t_h^\pm$ are determined by local (in each subdomain) behaviour
of the impedance.

\section{Coercivity estimate}

The coercivity of the skeleton equation \eqref{SkeletonEquation} is a
valuable feature because it garantees convergence of linear solvers.
To properly estimate the speed of convergence though, we need to bound this
coercivity estimate which is the focus of the present section.
We first establish an intermediate estimation. We will consider
\begin{equation*}
  \begin{aligned}
    & \mbVh(\Sigma)\times\mbVh(\Sigma)^*\;\text{equipped with}\\
    & \Vert (\ctrv,\ctrq)\Vert_{\Ts\times\Ts^{-1}}^2:=
    \Vert \ctrv\Vert_{\Ts}^2 + \Vert \ctrq\Vert_{\Ts^{-1}}^2
  \end{aligned}
\end{equation*}
which is the natural cartesian product norm. For this space, we have previously considered two
important subspaces: $\mathscr{X}_h(\Sigma) := \mbX_h(\Sigma)\times\mbX_h(\Sigma)^\circ$ in Lemma
\ref{CaracTransCond}, and the space of Cauchy data $\mathscr{C}_h(\Sigma)$ in Lemma \ref{InterpretationScattering}.
Besides $t_h^\pm$ equivalently defined by \eqref{BoundsTraceOperator} or \eqref{BoundsTraceOperator2}, 
we shall also rely on the inf-sup constant $\alpha_h$ and the continuity modulus $\Vert a\Vert$
defined by \eqref{MatrixInfSup}. 

\begin{prop}\label{EstimationProjectorCauchy}\quad\\
  Let Assumptions \ref{Hypo1}, \ref{Hypo2} and \ref{CoercivityImpedance} hold and 
  set $\mathscr{X}_h(\Sigma) := \mbX_h(\Sigma)\times\mbX_h(\Sigma)^\circ$.
  Then we have the following (a priori not orthogonal) direct sum
  \begin{equation*}
    \mbVh(\Sigma)\times\mbVh(\Sigma)^* =  \mathscr{X}_h(\Sigma)\oplus \mathscr{C}_h(\Sigma).
  \end{equation*}
  Moreover if $\mP: \mbVh(\Sigma)\times\mbVh(\Sigma)^*\to
  \mbVh(\Sigma)\times\mbVh(\Sigma)^*$ refers to the projector with
  $\mrm{Ker}(\mP) = \mathscr{X}_h(\Sigma)$ and $\mrm{Range}(\mP) = \mathscr{C}_h(\Sigma)$,
  we have
  \begin{equation*}
    \Vert \mP\Vert_{\Ts\times\Ts^{-1}}:=
    \sup_{(\ctrv,\ctrq)\in \mbVh(\Sigma)\times\mbVh(\Sigma)^*\setminus\{(0,0)\}}
    \frac{\Vert \mP(\ctrv,\ctrq)\Vert_{\Ts\times\Ts^{-1}}}{
      \Vert (\ctrv,\ctrq)\Vert_{\Ts\times\Ts^{-1}}}\leq \frac{(t_h^+)^2 + (2\Vert a\Vert/t_h^-)^2}{\alpha_h}
  \end{equation*}
\end{prop}
\noindent \textbf{Proof:}

First assume that $(\ctrp,\ctrq)\in \mathscr{X}_h(\Sigma)\cap \mathscr{C}_h(\Sigma)$. By
the very definition of $\mathscr{C}_h(\Sigma)$, there exists $\ctru\in\mbVh(\Omega)$ such that
$\Bh\ctru = \ctrp$ and $\Ah\ctru = \Bh^*\ctrq$. Since $\ctrp\in \mbX_h(\Sigma)$, we conclude that
$\ctru\in \mbX_h(\Omega)$. Next for any $\ctrv\in \mbX_h(\Omega)$ we have $\Bh\ctrv\in \mbX_h(\Sigma)$
and hence $\langle \Bh^*\ctrq,\ctrv\rangle = \langle \ctrq,\Bh\ctrv\rangle = 0$ as $\ctrq\in \mbX_h(\Sigma)^\circ$.
To sum up, we have $\ctru\in \mbX_h(\Omega)$ and $\langle \Ah\ctru,\ctrv\rangle = \langle\Bh^*\ctrq,\ctrv\rangle = 0$
for all $\ctrv\in\mbX_h(\Omega)$. Using \eqref{MatrixInfSup}, we deduce that $\ctru = 0$ hence $\ctrp =
\Bh\ctru = 0$ and $\Bh^*\ctrq = \Ah\ctru = 0\Rightarrow \ctrq = 0$ since $\Bh^*$ is one-to-one (because
$\Bh$ is onto). We have just established that
\begin{equation*}
  \mathscr{X}_h(\Sigma)\cap \mathscr{C}_h(\Sigma) = \{0\}.
\end{equation*}
Now we show that $\mathscr{X}_h(\Sigma) + \mathscr{C}_h(\Sigma) = \mbVh(\Sigma)\times\mbVh(\Sigma)^*$.
Pick an arbitrary pair $(\ctrp_{\dir},\ctrp_\neu)\in \mbVh(\Sigma)\times\mbVh(\Sigma)^*$.
Define $\tilde{\ctru}$ as the unique element of $\mbX_h(\Omega)$ such that
$\langle \Ah(\tilde{\ctru}+\Bh^\dagger\ctrp_{\dir}),\ctrv\rangle = \langle \ctrp_{\neu},\Bh(\ctrv)\rangle$
for all $\ctrv\in \mbX_h(\Omega)$. Taking account of \eqref{MatrixInfSup} then yields
\begin{equation}\label{Estimation1}
  \alpha_h\Vert \tilde{\ctru}\Vert_{\mbH^1(\Omega)}\leq
  (\Vert a \Vert/t_h^-)\Vert \ctrp_\dir\Vert_{\Ts}+t_h^+\Vert \ctrp_\neu\Vert_{\Ts^{-1}}
\end{equation}
Next let us set $\ctru = \tilde{\ctru}+\Bh^\dagger(\ctrp_\dir)$ and
$\ctru_\dir= \Bh(\ctru)= \Bh(\tilde{\ctru})+\ctrp_\dir$. For all
$\ctrv\in \mbVh(\Omega)$ satisfying $\Bh(\ctrv) = 0$ we have
$\langle \Ah\ctru,\ctrv\rangle = \langle \ctrp_{\neu},\Bh(\ctrv)\rangle = 0$ hence,
applying Lemma \ref{LiftingBoundaryFunctional} yields the existence of $\ctru_\neu$
such that $\langle \Ah\ctru,\ctrv\rangle = \langle \ctru_\neu,\Bh(\ctrv)\rangle\;
\forall \ctrv\in \mbVh(\Omega)$ which rewrites $\Ah\ctru = \Bh^*\ctru_\neu$.
In particular we have $\langle \ctru_\neu,\ctrq\rangle = \langle \Ah\ctru,\Bh^\dagger(\ctrq)\rangle$
for all $\ctrq\in \mbVh(\Sigma)$. From this we deduce the estimates
\begin{equation}\label{Estimation2}
  \begin{aligned}
    & \Vert \ctru\Vert_{\mbH^1(\Omega)} \leq \Vert \tilde{\ctru}\Vert_{\mbH^{1}(\Omega)} + (1/t_h^-) \Vert \ctrp_\dir\Vert_{\Ts} \\
    & \Vert \ctru_\dir\Vert_{\Ts}     \leq t_h^+\Vert \ctru\Vert_{\mbH^{1}(\Omega)}\\
    & \Vert \ctru_\neu\Vert_{\Ts^{-1}}  \leq (\Vert a\Vert/t_h^-)\Vert \ctru\Vert_{\mbH^1(\Omega)}
  \end{aligned}
\end{equation}
Now observe that, by construction, we have $(\ctru_\dir,\ctru_\neu)\in \mathscr{C}_h(\Sigma)$
see Proposition \ref{InterpretationScattering}. Besides $\ctrp_\dir-\ctru_\dir =
-\Bh(\tilde{\ctru})\in \mbX_h(\Sigma)$ since $\tilde{\ctru}\in \mbX_h(\Omega)$. On the other
hand $\langle \ctrp_\neu-\ctru_\neu,\Bh(\ctrv)\rangle = \langle \Ah(\ctru),\ctrv\rangle -
\langle \Ah(\ctru),\ctrv\rangle = 0$ for all $\ctrv\in \mbX_h(\Omega)$ and, since
$\mbX_h(\Sigma) = \Bh(\mbX_h(\Omega))$, we conclude that $\ctrp_\neu-\ctru_\neu\in\mbX_h(\Sigma)^\circ$.
In conclusion we have proved that $(\ctru_\dir,\ctru_\neu) - (\ctrp_\dir,\ctrp_\neu)\in \mathscr{X}_h(\Sigma)$
hence $(\ctru_\dir,\ctru_\neu) = \mP(\ctrp_\dir,\ctrp_\neu)$. To conclude the proof, there only remains to combine 
\eqref{Estimation1} and \eqref{Estimation2}.  \hfill $\Box$

\quad\\
The projection $\mP$ introduced in the previous result is intimately connected to the
inverse of the operator $\Id+\Pi\Sh$ of the skeleton formulation \eqref{SkeletonEquation}.
This is made apparent by the factorized form provided by the next theorem. 

\begin{thm}\label{FactorizationForm}\quad\\
  Let Assumptions \ref{Hypo1}, \ref{Hypo2} and \ref{CoercivityImpedance} hold.
  Define $\calT:\mbVh(\Sigma)\to \mbVh(\Sigma)\times\mbVh(\Sigma)^*$ by
  $\calT(\ctrv)=(\ctrv,-i\mT^*\ctrv)$ and
  $\calT':\mbVh(\Sigma)\times\mbVh(\Sigma)^*\to \mbVh(\Sigma)^*$ by 
  $\calT'(\ctrv,\ctrp)=\ctrp-i\mT(\ctrv)$.  Then we have  
  \begin{equation*}
    (\Id+\Pi\Sh)^{-1} = i\calT'\cdot\mP\cdot\calT\cdot\Ts^{-1}/2.
  \end{equation*}
\end{thm}
\noindent \textbf{Proof:}

Pick an arbitrary $\ctrf\in \mbVh(\Sigma)^*$ and define $\ctrp_{\dir} = \Ts^{-1}(\ctrf)/2$
and $\ctrp_\neu = -i\mT^*\Ts^{-1}(\ctrf)/2$, which simply rewrites $(\ctrp_\dir,\ctrp_\neu) =
\calT\cdot\Ts^{-1}(\ctrf)/2$. Next define $(\ctru_\dir,\ctru_\neu) = \mP(\ctrp_\dir,\ctrp_\neu)$
hence $(\ctru_\dir-\ctrp_\dir, \ctru_\neu-\ctrp_\neu)\in \mathscr{X}_h(\Sigma) = \mbX_h(\Sigma)\times
\mbX_h(\Sigma)^\circ$. On the other hand $(\ctru_\dir,\ctru_\neu)\in \mathscr{C}_h(\Sigma)$, so 
applying Lemma \ref{CaracTransCond} and Proposition \ref{InterpretationScattering} yields 
\begin{equation}\label{CalculInverse}
  \begin{aligned}
    & -(\ctru_\neu-\ctrp_\neu) + i\mT(\ctru_\dir-\ctrp_\dir) = \Pi(\;(\ctru_\neu-\ctrp_\neu) +i\mT^*(\ctru_\dir-\ctrp_\dir)\;)\\
    \iff\quad 
    & \ctru_\neu- i\mT\ctru_\dir + \Pi(\ctru_\neu+ i\mT^*\ctru_\dir) = (\Id + \Pi)\ctrp_\neu -i(\mT-\Pi\mT^*)\ctrp_\dir\\
    \iff\quad
    & (\Id+\Pi\Sh)(\ctru_\neu- i\mT\ctru_\dir) = (\Id + \Pi)\ctrp_\neu -i(\mT-\Pi\mT^*)\ctrp_\dir\\
    \iff\quad
    & (\Id+\Pi\Sh)\calT'\mP\calT\Ts^{-1}(\ctrf)/2 = (\Id + \Pi)\ctrp_\neu -i(\mT-\Pi\mT^*)\ctrp_\dir\\
  \end{aligned}
\end{equation}
Coming back to the definition of the exchange operator $\Pi$ given by Lemma \ref{ObliqueIsometry}, and 
setting $\mrm{Q} = \mR(\mR^*\mT^*\mR)^{-1}\mR^*$, we have $\Id+\Pi = 2\Ts\mrm{Q}$ and $\mT-\Pi\mT^*
= 2\Ts(\Id-\mQ\mT^*)$. Combining these identities with the definition of $(\ctrp_\dir,\ctrp_\neu)$ leads to
$(\Id + \Pi)\ctrp_\neu -i(\mT-\Pi\mT^*)\ctrp_\dir = -i\ctrf$. We have thus established the desired result.
\hfill $\Box$

\quad\\
Note that $\mathcal{T}'\cdot\mathcal{T} = 2\Ts/i$, so that $(i/2)\mathcal{T}\Ts^{-1}\mathcal{T}'$ is 
(a priori oblique) a projector, although this observation seems of no use in the present context.
We will exploit Theorem \ref{FactorizationForm} to establish a  coercivity estimate for
$\Id + \Pi\Sh$. The non-self adjoint part of the impedance operator $\mT$ shall arise naturally
in this analysis so we introduce a notation for bounding it
\begin{equation}\label{NonSelfAdjointEstimate}
  \begin{aligned}
    t_h^*
    & :=\sup_{\ctrv\in \mbVh(\Sigma)\setminus\{0\}}\Vert(\frac{\mT-\mT^*}{2})\ctrv\Vert_{\Ts^{-1}}/\Vert \ctrv\Vert_{\Ts}\\
    &  =\sup_{\ctrv\in \mbVh(\Sigma)\setminus\{0\}}\Big\vert \frac{\langle (\mT-\mT^*)\ctrv,\overline{\ctrv}\rangle}{
      \langle (\mT+\mT^*)\ctrv,\overline{\ctrv}\rangle}\Big\vert.
  \end{aligned}
\end{equation}

\begin{thm}\label{ExplicitLowerBound}\quad\\
  Let Assumptions \ref{Hypo1}, \ref{Hypo2} and \ref{CoercivityImpedance} hold.
  Then the inf-sup constant $\gamma_h$ from \eqref{CoercivityProperty} admits the following lower bound
  \begin{equation*}
    \gamma_h\geq \frac{2\alpha_h}{\lbr 1+(1+t_h^*)^2\rbr\,\lbr (t_h^+)^2+(2\Vert a\Vert/t_h^-)^2\rbr }
  \end{equation*}
\end{thm}
\noindent \textbf{Proof:}

To bound $\gamma_h$ from below, it suffices to bound $1/\gamma_h$ from above. On the other hand,
close inspection of the definition of $\gamma_h$ from \eqref{CoercivityProperty} shows that
$1/\gamma_h$ is the continuity modulus of $(\Id-\Pi\Sh)^{-1}$ which writes
\begin{equation}\label{NormInverseSkeletonOperator}
  \frac{1}{\gamma_h} = \sup_{\ctrp\in \mbVh(\Sigma)^*\setminus\{0\}}
  \frac{\Vert (\Id + \Pi\Sh)^{-1}\ctrp\Vert_{\Ts^{-1}}}{\Vert\ctrp\Vert_{\Ts^{-1}}}.
\end{equation}
Next we need to obtain upper bounds for the continuity modulus of $\mathcal{T},\mathcal{T}'$
and $\mP$. A direct estimation using \eqref{NonSelfAdjointEstimate} yields
\begin{equation}\label{EstimationContinuiteOperatorT}
  \begin{aligned}
    \Vert \mathcal{T}'(\ctrv,\ctrq)\Vert_{\Ts^{-1}}
    & \leq \Vert \ctrq\Vert_{\Ts^{-1}} + (1+t_h^*)\Vert \ctrv\Vert_{\Ts} \\
    & \leq (1+(1+t_h^*)^2)^{1/2}\Vert (\ctrv,\ctrq)\Vert_{\Ts\times\Ts^{-1}}\\
    \Vert \mathcal{T}(\ctrv)\Vert_{\Ts\times\Ts^{-1}}
    & \leq (\Vert \ctrv\Vert_{\Ts}^2 + \Vert \mT^*(\ctrv)\Vert_{\Ts^{-1}}^2)^{1/2}\\
    & \leq (1+(1+t_h^*)^2)^{1/2}\Vert \ctrv\Vert_{\Ts}
  \end{aligned}
\end{equation}
In the estimate above we used the elementary inequalities
$\Vert \mT(\ctrv)\Vert_{\Ts^{-1}} \leq (1+t_h^*)\Vert \ctrv\Vert_{\Ts}$
and similarly for $\mT^*$, and  $\Vert \Ts(\ctrv)\Vert_{\Ts^{-1}} =
\Vert \ctrv\Vert_{\Ts}$. Plugging the factorized form of Theorem
\ref{FactorizationForm}, combined with \eqref{EstimationContinuiteOperatorT}
and Proposition \ref{EstimationProjectorCauchy}, into
\eqref{NormInverseSkeletonOperator} yields the desired estimate.
\hfill $\Box$

\quad\\
The next result, that agrees with the estimate provided by
\cite[Prop. 10.4]{claeys2020robust}, yields another variant of this
coercivity bound that appears sharper in certain cases. 

\begin{prop}\quad\\
  Let Assumptions \ref{Hypo1}, \ref{Hypo2} and \ref{CoercivityImpedance} hold,
  and suppose in addition that $\mT = \mT^* = \Ts$. Then we have the following
  estimate
  \begin{equation}
    \gamma_h = \inf_{\ctrq\in \mbVh(\Sigma)^*\setminus\{0\}}
    \frac{\Vert (\Id +\Pi\Sh)\ctrq\Vert_{\Ts^{-1}}}{\Vert \ctrq\Vert_{\Ts^{-1}}}
    \geq \frac{1}{\Vert \mP\Vert_{\Ts\times\Ts^{-1}}}.
  \end{equation}
\end{prop}
\noindent \textbf{Proof:}

Since $\mT = \Ts$ we have $\Vert \mathcal{T}'(\ctrv,\ctrp)\Vert_{\Ts^{-1}}\leq \sqrt{2}
\Vert (\ctrv,\ctrp)\Vert_{\Ts\times\Ts^{-1}}$ for all $\ctrv\in\mbVh(\Sigma),\ctrp\in\mbVh(\Sigma)^*$.
Moreover a direct calculation yields $\Vert \mathcal{T}\Ts^{-1}(\ctrp)\Vert_{\Ts\times\Ts^{-1}}^2 =
2\Vert\ctrp\Vert_{\Ts^{-1}}^2$. Using the definition of the
continuity modulus of $\mP$, as well as Theorem \ref{FactorizationForm}
and the definition of $\gamma_h$, we obtain the inequality
\begin{equation*}
  \begin{aligned}
    \frac{1}{\gamma_h} = \sup_{\ctrp\in \mbVh(\Sigma)^*\setminus\{0\}}
    \frac{\Vert \mathcal{T}'\mP\mathcal{T}\Ts^{-1}(\ctrp)\Vert_{\Ts\times\Ts^{-1}}
    }{2 \Vert \ctrp\Vert_{\Ts^{-1}}}
     & \leq
    \sup_{\ctrp\in \mbVh(\Sigma)^*\setminus\{0\}}
    \frac{\Vert \mP\mathcal{T}\Ts^{-1}(\ctrp)\Vert_{\Ts\times\Ts^{-1}}
    }{\Vert \mathcal{T}\Ts^{-1}(\ctrp)\Vert_{\Ts\times\Ts^{-1}}}\\
    & \leq \Vert \mP\Vert_{\Ts\times\Ts^{-1}}.
  \end{aligned}
\end{equation*}

\hfill $\Box$

\section{Purely local exchange operator}\label{PurelyLocalImpedance}

In this section, we wish to show that the theory of the present contribution
covers other pre-existing variants of OSM that involve non-self-adjoint
impedance operators. We will in particular establish connections
with the DDM strategies described in \cite{zbMATH07197799} and
\cite{despres:hal-03230250}. Additionally let us recall that,
as explained in \cite[Sect.12]{claeys2020robust}, our theoretical framework
also covers the initial Despr\'es algorithm of
\cite{MR1291197,MR1105979,MR1071633,MR1227838}.

\quad\\
Applying the exchange operator introduced in Section
\ref{ReformulationTransmissionConditions} (i.e.  the operation
$\ctrp\mapsto \Pi(\ctrp)$) is computationally non-trivial
because the expression of $\Pi$ involves the inverse operator $(\mR^*\mT^*\mR)^{-1}$
which is a priori not known explicitely. However in the particular case where
the impedance is diagonal, calculations related to the exchange operator
become much simpler. The present section will focus on the even simpler situation
where the impedance is associated to the identity matrix which
is already an interesting and instructive special case.

\subsection{Derivation of the swapping operator}
Let us denote $\dof(\Omega)$ the set of degrees of freedom of the $\mathbb{P}_k$-Lagrange
discretization over the whole computational domain, which will be considered a point cloud, and set
$\dof(\Gamma_j) =  \Gamma_j\cap \dof(\Omega)$ and $\dof(\Sigma) =  \Sigma\cap \dof(\Omega)$.
In the present section we will assume that $\mT = \mD = \mrm{diag}_{j=1\dots \mJ}(\mD_j)$
where $\mD:\mbVh(\Sigma)\to \mbVh(\Sigma)^*$ is defined by 
\begin{equation}\label{DefDiagOperator}
  \begin{aligned}
    & \langle \mD(\ctru),\ctrv\rangle :=
    \langle \mD_1(u_1),v_1\rangle + \dots +
    \langle \mD_\mJ(u_\mJ),v_\mJ\rangle\\
    & \text{with}\quad \langle \mD_j(u_j),v_j\rangle := \sum_{\bx\in \dof(\Gamma_j)} u_j(\bx)v_j(\bx)
  \end{aligned}
\end{equation}
for any $\ctru,\ctrv\in \mbVh(\Sigma)$ with $\ctru = (u_1,\dots, u_\mJ)$ and
$\ctrv = (v_1,\dots, v_\mJ)$. Obviously this is a self-adjoint operator $\mD =\mD^*$.
Algebraically, each $\mD_j$ simply corresponds to an identity matrix. We shall denote
$\Pi_{\loc}$ the associated exhange operator which is defined by $\Pi_{\loc} = 2\mD(\mR^*\mD\mR)^{-1}\mR^*-\Id$
in accordance with Lemma \ref{ObliqueIsometry}. 

\quad\\
To obtain a more explicit expression for the exchange operator $\Pi_\loc$, we start from an arbitrary
$\ctrp = (p_1,\dots,p_\mJ)\in\mbVh(\Sigma)^*$ and we decompose each local component
$p_j\in\Vh(\Gamma_j)^*,j=1\dots \mJ$ according to its coordinates in the canonical dual basis
which writes
\begin{equation*}
  \langle p_j,v\rangle = \sum_{\bx\in \dof(\Gamma_j)} p_j(\bx)v(\bx)\quad \forall v\in \Vh(\Gamma_j).
\end{equation*}
In other words, defining $\delta_{\bx}\in \Vh(\Gamma_j)^*$ by $\langle \delta_{\bx},v\rangle := v(\bx)$,
we have $p_j = \sum_{\bx\in \dof(\Gamma_j)}p_j(\bx)\delta_{\bx}$. We derive an expression of $\Pi(\ctrp)$ in
terms of the coefficients $p_j(\bx)$.  Examining the action of the operators $\mR^*:\mbVh(\Sigma)^*\to
\Vh(\Sigma)^*$ and $\mR^*\mD\mR:\mbVh(\Sigma)\to\mbVh(\Sigma)^*$ we obtain the following:
for any $u,v\in \Vh(\Sigma)$
\begin{equation}
 \begin{aligned}
    \langle \mR^*(\ctrp),v\rangle  =   \langle \ctrp,\mR(v)\rangle
    = \sum_{\bx\in \dof(\Sigma)}v(\bx)\sum_{j:\bx\in\Gamma_j} p_j(\bx)\\
    \langle \mD\mR(u),\mR(v)\rangle =
    \sum_{\bx\in \dof(\Sigma)} \#\{j:\bx\in\Gamma_j\} \;u(\bx)v(\bx)
  \end{aligned}
\end{equation}
In particular $\mR^*\mD\mR$ appears diagonal. Now plugging the previous
expressions into the definition $\Pi_\loc := 2\mD\mR(\mR^*\mD\mR)^{-1}\mR^*-\Id$,
for any $\ctrv = (v_1,\dots,v_\mJ)\in\mbVh(\Sigma)$ we obtain
\begin{equation}\label{SwappingOperator}
  \begin{aligned}
    & \langle \Pi_\loc(\ctrp),\ctrv\rangle =\\
    & \sum_{j=1\dots \mJ}\sum_{\bx\in\dof(\Gamma_j)}v_j(\bx)
    \Big\lbr -p_j(\bx) + \frac{2}{\#\{n:\bx\in\Gamma_n\}}\sum_{k:\bx\in\Gamma_k}p_k(\bx)\Big\rbr
  \end{aligned}
\end{equation}
Let us inspect this expression, assuming for a moment that the subdomain partition
does not contain any  cross-point. In this case we have
$\#\{n:\bx\in\Gamma_n\} = 2\;\forall \bx\in\dof(\Sigma)$ and the operator $-\Pi_\loc$
simply consists in swapping the unknowns from both
sides of each interface. This operator $\Pi_\loc$ should thus be understood
as the exchange operator that is used in the rest of the literature on
Optimized Schwarz Methods. This is indeed the operator considered in
\cite{zbMATH05161203,despres:hal-03230250,despres:hal-02612368,MR3450068,MR1764190,
  lecouvez:tel-01444540,MR3989867,MR3549886,MR2097757,MR1924414,MR1227838,MR1291197,
  MR1071633,zbMATH07197799,Boubendir2012,ElBouajaji2015,Modave2020b,Bendali2006}
for enforcing coupling between neighbouring subdomains. Expression
\eqref{SwappingOperator} can be re-arranged so as to take a more symetric form
\begin{equation}\label{SwappingOperator2}
  \langle \Pi_\loc(\ctrp),\ctrv\rangle =  -\langle \ctrp,\ctrv\rangle
  + 2\sum_{\bx\in\dof(\Sigma)}\frac{\Big\lbr\sum_{j:\bx\in\Gamma_j}p_j(\bx)\Big\rbr
  \Big\lbr\sum_{k:\bx\in\Gamma_k}v_k(\bx)\Big\rbr}{\#\{n:\bx\in\Gamma_n\}}.
\end{equation}
Applying Lemma \ref{ProprieteEchangeOrthogonal} with $\mT = \mD$ yields that $\Pi_\loc^*\mR = \mR$, $\Pi_\loc^2 = \Id$ i.e.
swapping traces twice leaves them unchanged, and $\Pi_\loc(\ctrq) = -\ctrq\;\forall q\in\mbX_h(\Sigma)^\circ$.

\begin{remark}
 Consider the case where the impedance is chosen as $\mT = z\mD$ where $z\in\CC$ and
  $\mD$ is the diagonal operator defined by \eqref{DefDiagOperator}. Define $\mP_\loc:=(\Id+\Pi_{\loc})/2
  = \mD\mR(\mR^*\mD\mR)^{-1}\mR^*$. Following Remark \ref{NonInvolutive}, we see that the exchange
  operator associated to this choice of impedance is given by
  \begin{equation*}
    \Pi = \Big(\frac{z}{\vert z\vert}\Big)^{2}\mP_{\loc} - (\Id - \mP_{\loc})
  \end{equation*}
  It is purely local an exchange operator (i.e. it only couples degrees of freedom that geometrically coincide),
  and at the same time $\Pi\neq \Pi_{\loc}$ unless $z = 1$
  (under Assumption \ref{CoercivityImpedance}). Interestingly, if $z = \vert z\vert \exp(im\pi/(2p+1))$
  for $m,p\in\mathbb{N}$, then $(z/\vert z\vert)^{2(2p+1)} = 1, (-1)^{2p+1} = -1$ hence, in this case, we
  have the identity $\Pi^{2p+1}=\Pi_{\loc}=\Pi_{\loc}^{2p+1}$.
\end{remark}

\subsection{A criterion for locality}
The expression \eqref{SwappingOperator} and \eqref{SwappingOperator2} are fully
explicit so that applying $\Pi_\loc$ is fast and straightforward. We proved that
the exchange operator induced by $\mD$ is  $\Pi_\loc = 
2\mD\mR(\mR^*\mD\mR)^{-1}\mR^*-\Id$. But other choices of impedance can induce the same
exchange operator. A natural question arises then to determine those impedance
operators that lead to $\Pi_\loc$ as exchange operator. The next lemma provides
an explicit criterion for this.

\begin{lem}\label{LocalityCriterionExchangeOperator}\quad\\
  For any linear map $\mT:\mbVh(\Sigma)\to \mbVh(\Sigma)^*$
  satisfying Assumption \ref{CoercivityImpedance} we have
  \begin{equation*}
    \Pi_\loc\mT^*\mR = \mT\mR  \iff  \Pi_\loc = (\mT+\mT^*)\mR(\mR^*\mT^*\mR)^{-1}\mR^*- \Id.
    \end{equation*}
\end{lem}
\noindent\textbf{Proof:}

Assume first that $\Pi_\loc = (\mT+\mT^*)\mR(\mR^*\mT^*\mR)^{-1}\mR^*- \Id$. Then a
direct calculus shows that $\Pi_\loc\mT^*\mR = 2\Ts\mR(\mR^*\mT^*\mR)^{-1}(\mR^*\mT^*\mR)- \mT^*\mR
= (\mT+\mT^*)\mR - \mT^*\mR= \mT\mR$.

Reciprocally assume that $\Pi_\loc\mT^*\mR = \mT\mR$ holds. Pick an arbitrary $\ctrp\in \mbVh(\Sigma)^*$
and let $u\in \Vh(\Sigma)$ solve $\langle \mT^*\mR(u),\mR(v)\rangle =
\langle \ctrp,\mR(v)\rangle\;\forall v\in\Vh(\Sigma)$,
which is uniquely solvable a variational problem thanks to the coercivity of $\mT^*$ given by Assumption
\ref{CoercivityImpedance}. Then we have $u = (\mR^*\mT^*\mR)^{-1}\mR^*\ctrp$ and, by construction,
$\ctrq = \ctrp - \mT^*\mR(u)\in \mbX_h(\Sigma)^\circ$. From this we conclude
\begin{equation*}
  \begin{aligned}
    \Pi_\loc(\ctrp)
    & = \Pi_\loc(\ctrq) + \Pi_\loc\mT^*\mR(u)\\
    & = -\ctrq+\mT\mR(u) = -\ctrq -\mT^*\mR(u) + (\mT+\mT^*)\mR(u)\\
    & = -\ctrp + (\mT+\mT^*)\mR(\mR^*\mT^*\mR)^{-1}\mR^*\ctrp
  \end{aligned}
\end{equation*}
In the above calculus we have used the fact that $\Pi_\loc(\ctrq) = -\ctrq\;\forall\ctrq\in \mbX_h(\Sigma)^\circ$.
Since the tuple of traces $\ctrp$ was chosen arbitrarily in $\mbVh(\Sigma)^*$, we have proved the desired result.
\hfill $\Box$

\quad\\
As a corollary, the criterion exhibited in Lemma \ref{LocalityCriterionExchangeOperator} can
be simplified, taking the form of a commutation identity. This is not an equivalence anymore though.

\begin{cor}\label{CriterionLocalityCommutativity}\quad\\
  For any linear map $\mT:\mbVh(\Sigma)\to \mbVh(\Sigma)^*$
  satisfying Assumption \ref{CoercivityImpedance} we have
  \begin{equation*}
    \Pi_\loc\mT = \mT^*\Pi_\loc^*\quad \Longrightarrow\quad
    \Pi_\loc = (\mT+\mT^*)\mR(\mR^*\mT^*\mR)^{-1}\mR^*-\Id.
  \end{equation*}
\end{cor}
\noindent\textbf{Proof:}

Observe that $(\Pi_\loc)^2 = \Id$, hence multiplying  $\Pi_\loc\mT = \mT^*\Pi_\loc^*$
on the left by $\Pi_\loc$ and on the right by  $\Pi_\loc^*$ yields    $\mT\Pi_\loc^* = \Pi_\loc\mT^*$.
Since we have $\Pi_\loc^*\mR = \mR$ systematically according to Lemma \ref{ProprieteEchangeOrthogonal},
we conclude that $\Pi_\loc\mT^*\mR = \mT\Pi_\loc^*\mR = \mT\mR$. There only remains to apply Lemma
\ref{LocalityCriterionExchangeOperator} which yields the desired result. \hfill $\Box$

\begin{remark}
Like in Example \ref{ExampleEMDAOO0} assume a decomposition in two subdomains
  with no cross-point, see Fig \ref{MainFig} (b). In this situation, the impedance operator associated to
  OO0 and EMDA strategies both take the form $\mT = z\mT_{\textsc{r}}$ where
  $\mT_{\textsc{r}}$ stems from surface mass matrices on the $\Gamma_j$'s. Then we know
  from \cite[Sect.12]{claeys2020robust} that, when there is no cross point, we have the relation
  $\Pi_{\loc}\mT_{\textsc{r}} = \mT_{\textsc{r}}\Pi_{\loc}^*$ i.e. the exchange
  operator associated to $\mT_{\textsc{r}}$ is $\Pi_{\loc}$. Following Remark \ref{NonInvolutive},
  the definition
  \begin{equation}\label{RightExchangeOperator}
    \Pi = \Big(\frac{z}{\vert z\vert}\Big)^2\mP_{\loc} - (\Id - \mP_{\loc})
  \end{equation}
  should be taken for the exchange operator. In particular, in the case of OO0
  we have $z/\vert z\vert = \exp(i\pi/4)$, hence $(z/\vert z\vert)^2 = i$ and 
  $\Pi = i\mP_{\loc} - (\Id - \mP_{\loc})$. We see that $\Pi^{2} = -\Pi_{\loc}$
  and $\Pi \neq \Pi_{\loc}$. 
\end{remark}

\begin{remark}
Let us now reconsider the situation of Example \ref{LienLecouvezCollinoJoly},
    again with a geometrical configuration depicted in Figure \ref{MainFig} (b).
    Choosing an impedance as in \eqref{ImpedanceCollinoJolyLecouvez}, for
    $\bu = (u_1,u_2), \bv = (v_1,v_2)$ with $u_j,v_j\in\Vh(\Gamma_j)$,
    since $\mT_e = \mT_e^{*}$ we have
  \begin{equation*}
    \begin{aligned}
      \langle \mT^{*}(\bu),\bv\rangle
      & = \langle \mT_1^*(u_1),v_1\rangle +
      \langle \mT_1(u_2\vert_{\Gamma_1}),v_2\vert_{\Gamma_1}\rangle +
      \langle \mT_e(u_2\vert_{\Gamma_e}),v_2\vert_{\Gamma_e}\rangle \\
      \langle \Pi_{\loc}\mT(\bu),\bv\rangle
      & = \langle \mT_1(u_1),v_2\vert_{\Gamma_1}\rangle +
      \langle \mT_1^*(u_2\vert_{\Gamma_1}),v_1\rangle +
      \langle \mT_e(u_2\vert_{\Gamma_e}),v_2\vert_{\Gamma_e}\rangle \\
      \langle \mT^*\Pi_{\loc}^*(\bu),\bv\rangle
      & = \langle \mT_1^*(u_2\vert_{\Gamma_1}),v_1\rangle +
      \langle \mT_1(u_1),v_2\vert_{\Gamma_1}\rangle +
      \langle \mT_e(u_2\vert_{\Gamma_e}),v_2\vert_{\Gamma_e}\rangle
    \end{aligned}
  \end{equation*}
  In this situation, the criterion  $\Pi_\loc\mT = \mT^*\Pi_\loc^*$ is satisfied.
  The choice of impedance  \eqref{ImpedanceCollinoJolyLecouvez} is the one
  considered in \cite[\S 2.1.2]{zbMATH07197799} where the analysis is conducted 
  with $\Pi_{\loc}$ as exchange operator. The criterion $\Pi_\loc\mT = \mT^*\Pi_\loc^*$
  thus shows how our theory recovers the results of \cite{zbMATH07197799}.
\end{remark}

\begin{remark}
The criterion $\Pi_\loc\mT = \mT^*\Pi_\loc^*$
  also directly matches the compatibility assumption of \cite[Def.12 \& Lem.14]{despres:hal-03230250}
  where a particular choice of impedance fulfilling Assumption \ref{CoercivityImpedance} is considered,
  see \cite[Prop.4]{despres:hal-03230250}.
\end{remark}

\noindent 
The previous remarks show that our theory covers the strategies considered in \cite{zbMATH07197799}
and \cite{despres:hal-03230250} embedding them in a more general framework and, in passing, provides
a refined convergence estimate for them through application of Theorem \ref{ExplicitLowerBound}. 

\quad\\
In conclusion, we would like to point that, in the special case of a self-adjoint impedance,
the criterion provided by Corollary \ref{CriterionLocalityCommutativity} does give rise to
an equivalence.

\begin{cor}\quad\\
  For any linear map $\mT:\mbVh(\Sigma)\to \mbVh(\Sigma)^*$ 
  satisfying Assumption \ref{CoercivityImpedance}, and that is in addition self-adjoint
  $\mT = \mT^*$, we have
  \begin{equation*}
    \Pi_\loc\mT = \mT\Pi_\loc^*\quad \iff\quad
    \Pi_\loc = 2\mT\mR(\mR^*\mT\mR)^{-1}\mR^*- \Id.
  \end{equation*}
\end{cor}
\noindent\textbf{Proof:}

Direct calculus indicates clearly that, if $\Pi_\loc = 2\mT\mR(\mR^*\mT\mR)^{-1}\mR^*- \Id$,
then $\Pi_\loc\mT = \mT\Pi_\loc^*$. The reciprocal follows from Corollary \ref{CriterionLocalityCommutativity}.
\hfill $\Box$

\section{Numerical illustration}

\begin{figure}[!]
  \begin{center}
    \begin{subfigure}[b]{0.4\textwidth}
      \centering
      \includegraphics[height=4cm]{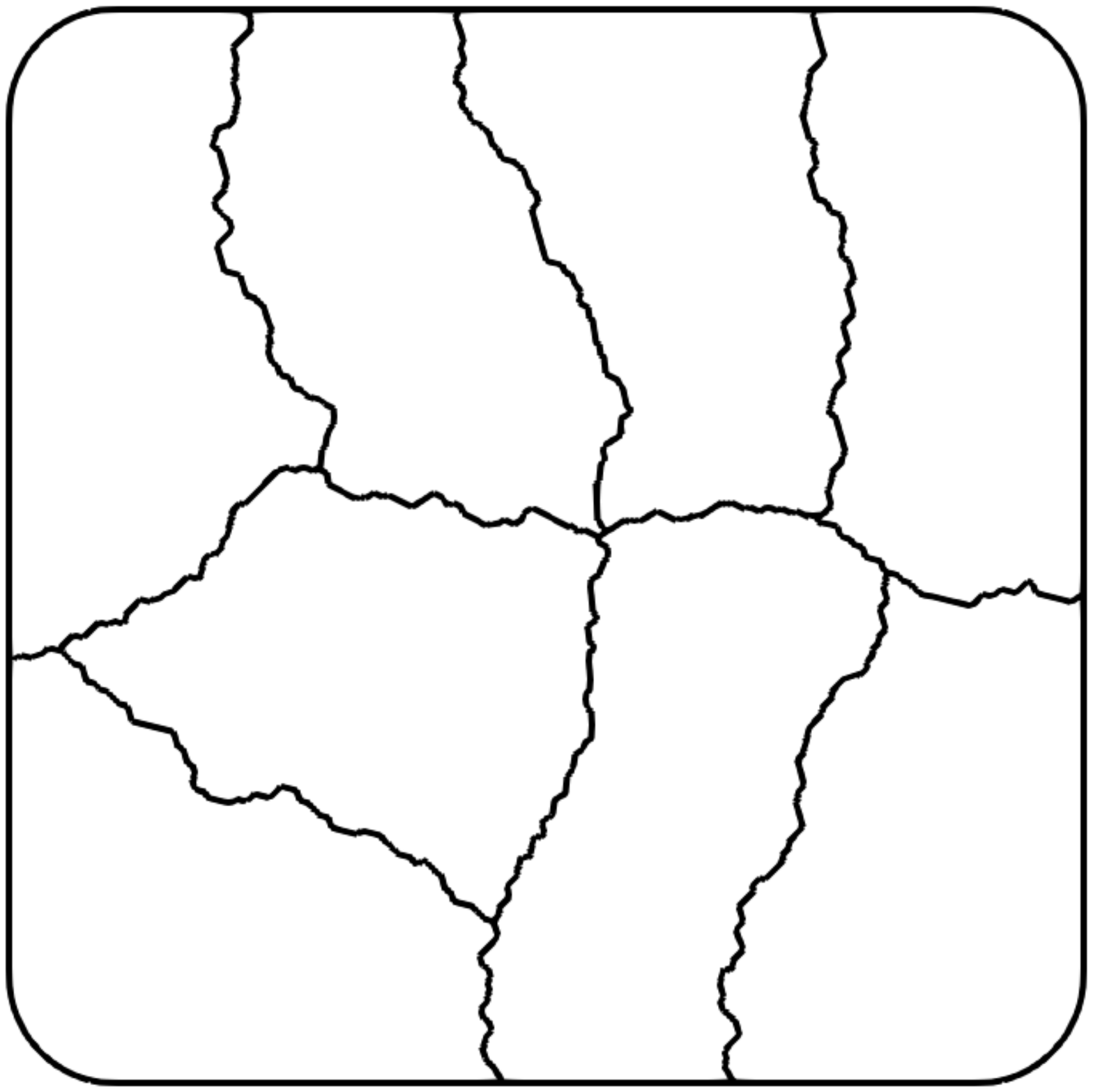}
      \caption{}
    \end{subfigure}
    \begin{subfigure}[b]{0.25\textwidth} 
      \centering
      \raisebox{0.4cm}{\includegraphics[height=3.25cm]{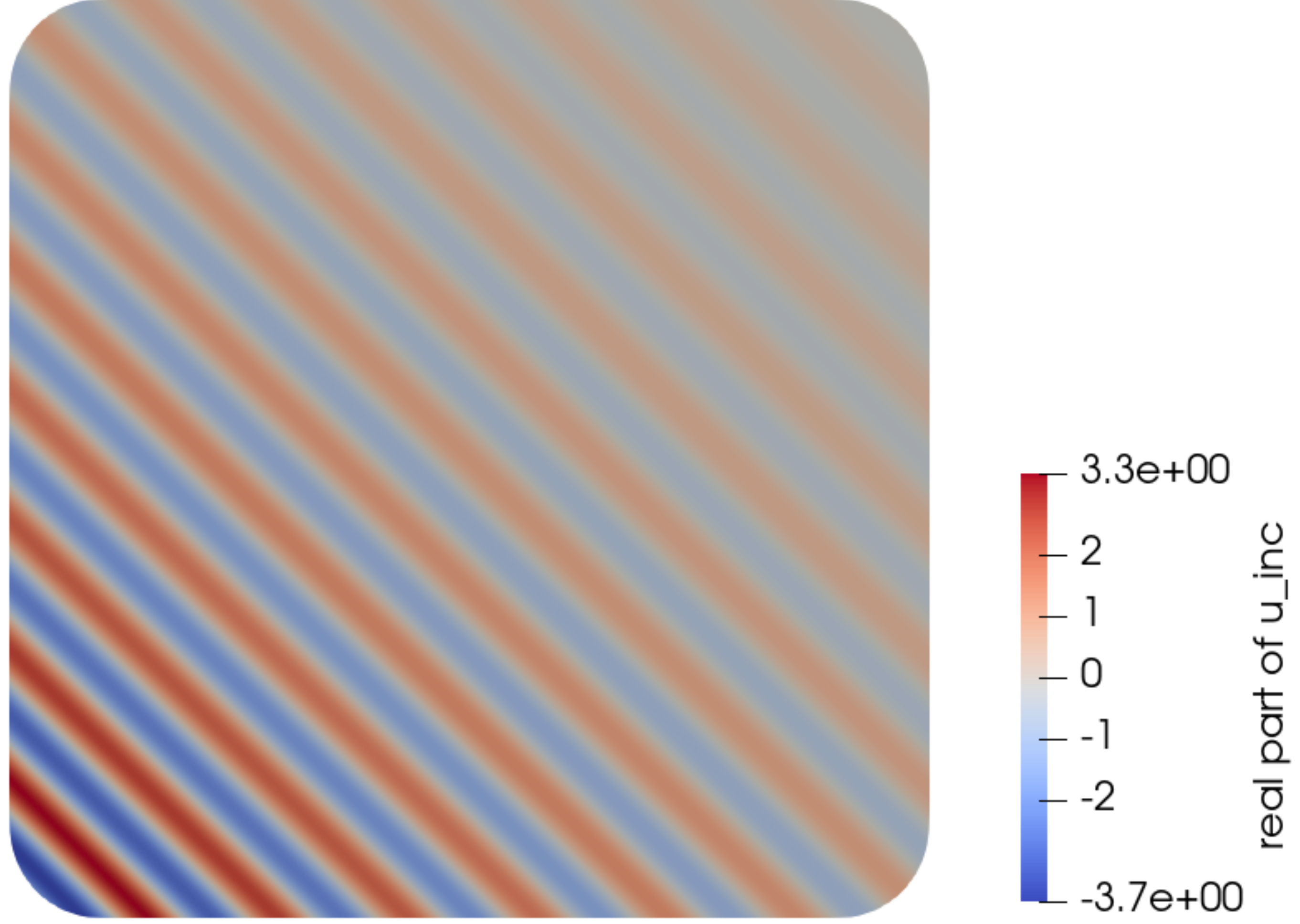}}
      \caption{}
    \end{subfigure}
    \caption{(a) non-overlapping partition, (b)  real part of the solution to be computed.}\label{ComputationalDomain}
  \end{center}
\end{figure}

This final section reports on a numerical example illustrating the theoretical
convergence results we have established previously. Our primal aim is to confirm the convergence
of a Richardson linear solver applied to an instance of the skeleton formulation \eqref{SkeletonEquation}.
In addition, we shall briefly examine whether a better convergence can be obtained by considering
non-HPD impedance operators. The geometry of the computational domain $\Omega$, a square of side length 2
centered at $\bx = 0$ with rounded corners, is depicted in Figure \ref{ComputationalDomain} together
with its partitioning.

\noindent 
  We target the numerical solution to the boundary value problem \eqref{InitialPb} in $\RR^d = \RR^2$
with $\mu = 1$ and a constant wave number $\kappa = 2\pi/\lambda + i \in\CC$ with $\lambda = 0.2$.
Regarding the source terms, we take $f = 0$ and $g = \partial_{\bn}u_{\mrm{inc}}$ where
$u_{\mrm{inc}}(\bx) = \exp(i\kappa \bd\cdot\bx)$ with $\bd = (1/\sqrt{2}, 1/\sqrt{2}, 0)$.
This boundary value problem is discretized with
$\mathbb{P}_1-$Lagrange finite elements on a triangular mesh with
23634 triangles and 12010 vertices generated by means of
\textsc{gmsh}\footnote{\texttt{https://gmsh.info/}}, and partitioned
in 8 subdomains using
\textsc{metis}\footnote{\texttt{https://github.com/KarypisLab/METIS}}.
We consider the skeleton formulation \eqref{SkeletonEquation} (with a right-hand
side stemming from $u_{\mrm{inc}}$) associated to three different choices of
impedance $\mT = \mrm{diag}(\mT_1,\dots,\mT_\mJ)$:\\[-5pt]
\begin{itemize}
\item \textbf{choice 1:} $\langle \mT_j(u),v\rangle = \kappa\int_{\Gamma_j}u v\,d\sigma$,
\item \textbf{choice 2:} $\langle \mT_j(u),v\rangle = \frac{1}{2\vert\kappa\vert}\int_{\Gamma_j}\nabla_{\Gamma_j}u \cdot\nabla_{\Gamma_j}v\,d\sigma
  + \vert\kappa\vert\int_{\Gamma_j}u v\,d\sigma$,
\item \textbf{choice 3:} $\langle \mT_j^{\theta}(u),v\rangle = \exp(-i\theta)\Big(\frac{1}{2\vert\kappa\vert}
  \int_{\Gamma_j}\nabla_{\Gamma_j}u \cdot\nabla_{\Gamma_j}v\,d\sigma
  + \vert\kappa\vert\int_{\Gamma_j}u v\,d\sigma\Big)$.  
\end{itemize}
In the expressions above $\nabla_{\Gamma_j}$ refers to the tangential gradient on $\Gamma_j$.
As opposed to Choice 2, Choices 1 and 3 are not HPD
impedance operators. The exchange operator $\Pi$ is implemented 
according to the formula of Lemma \ref{CaracTransCond}. This
formula involves the term $(\mR^*\mT^*\mR)^{-1}$ which requires solving
a linear system on the skeleton for each evaluation of the
matrix-vector product $\ctrp\mapsto \Pi(\ctrp)$. For our
implementation, this linear system is solved by means of
\textsc{umfpack}\footnote{\texttt{https://people.engr.tamu.edu/davis/suitesparse.html}}.

\begin{figure}
  \begin{minipage}{0.45\linewidth}
    \centering
    \includegraphics[height=5cm]{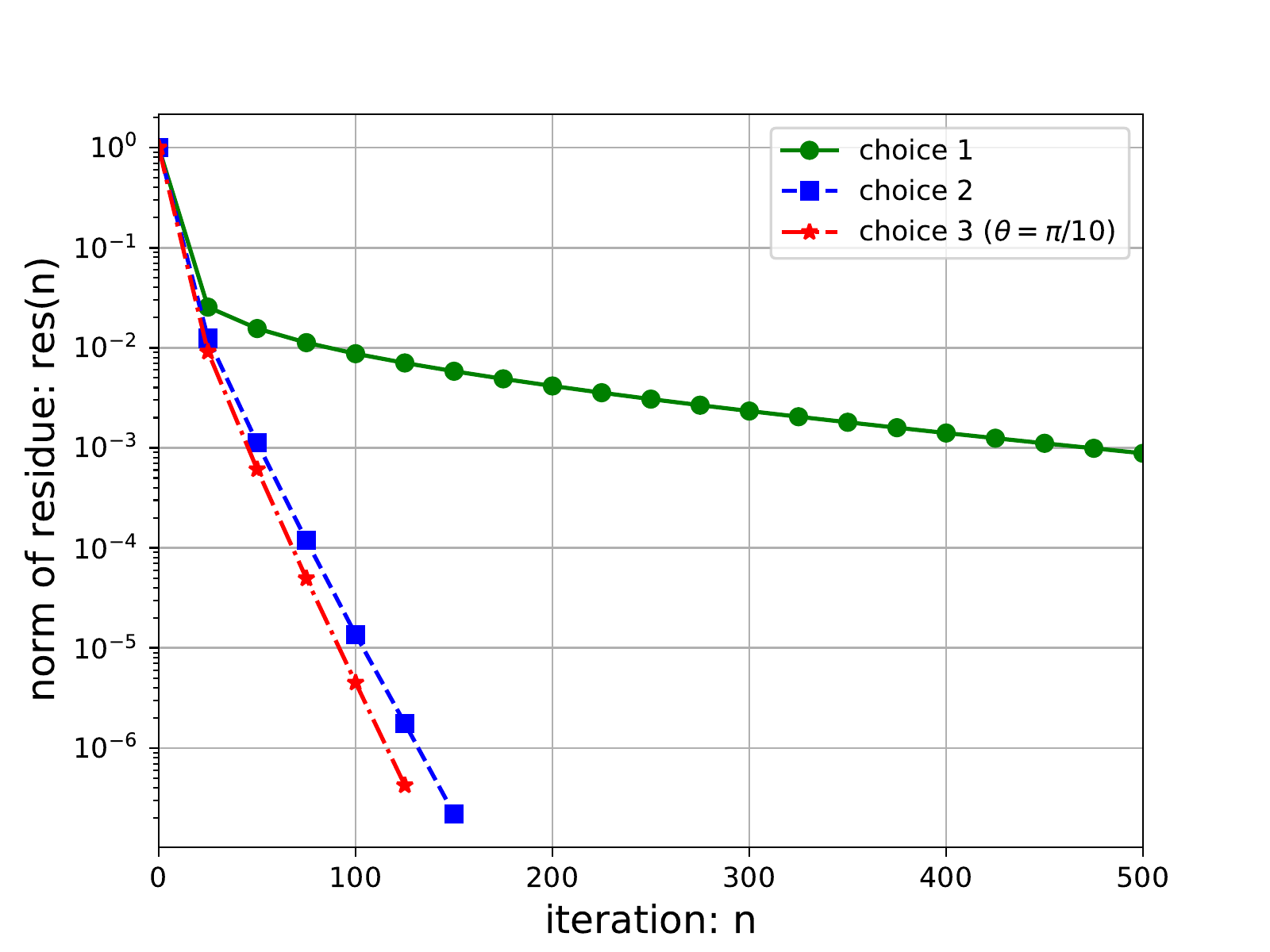}
    \caption{Norm of the residue vs iteration}\label{ConvergenceHistory}
  \end{minipage}
  \hspace{1cm}
  \begin{minipage}{0.4\linewidth}
    \centering
    \includegraphics[height=5cm]{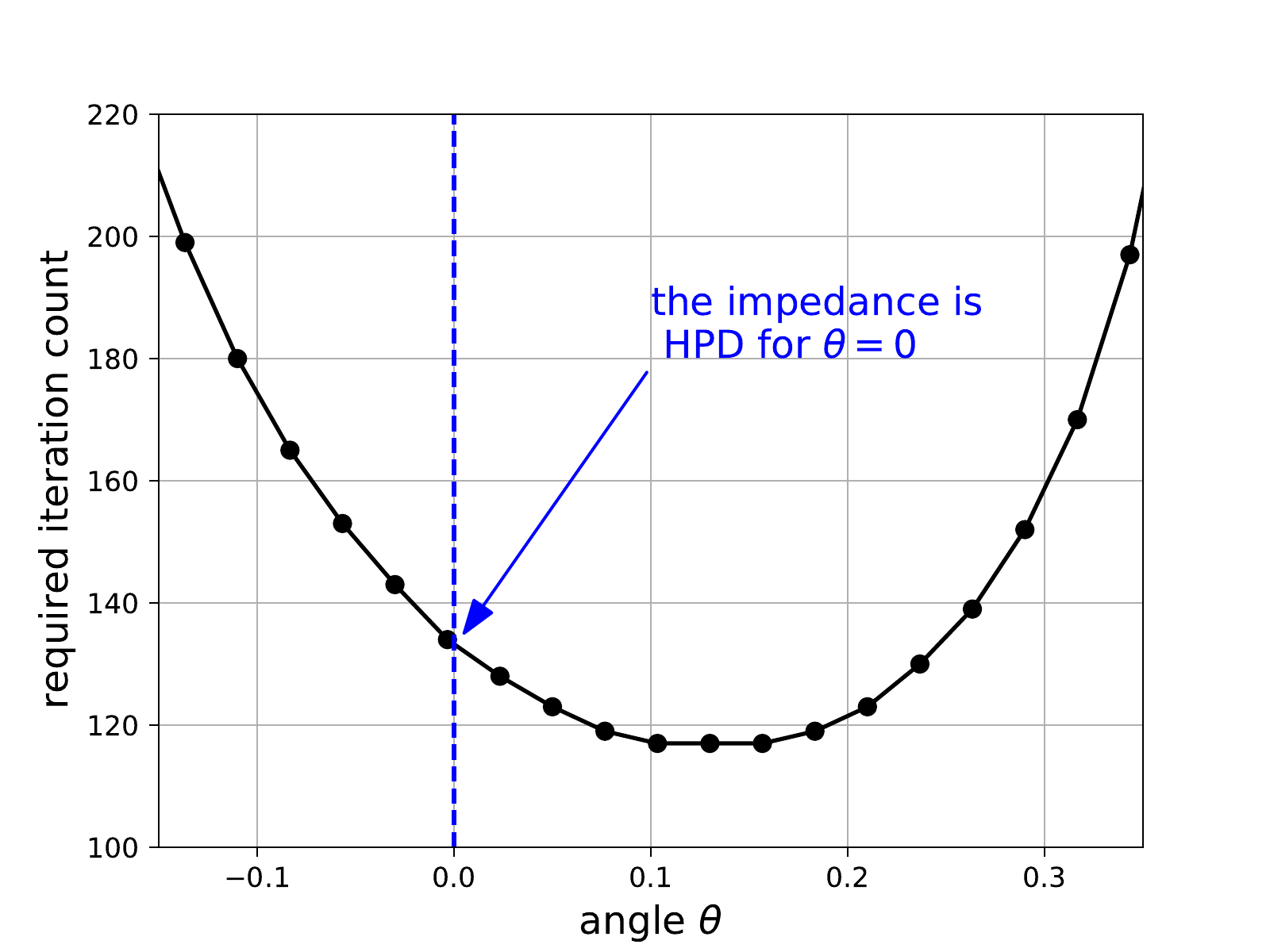}
    \caption{Iteration count vs angle $\theta$}\label{RequiredIter}
  \end{minipage}
\end{figure}

The overall skeleton formulation \eqref{SkeletonEquation} is solved
with a Richardson solver i.e. we compute the sequence of iterates
$\ctrq^{(n)}$ starting at $\ctrq^{(0)} = 0$ and then defined by
$\ctrq^{(n+1)} = \ctrq^{(n)} + \alpha(\ctrg-(\Id+\Pi\mS)\ctrq^{(n)})$
with relaxation parameter $\alpha = 1/\sqrt{2}$. The solver
is stopped when the following residual norm passes below $10^{-6}$:
\begin{equation*}
  \mrm{res}(n):=\Vert \ctrg-(\Id+\Pi\mS)\ctrq^{(n)}\Vert_{\Ts^{-1}}/\Vert \ctrg\Vert_{\Ts^{-1}}.
\end{equation*}
Figure \ref{ConvergenceHistory} confirms the systematic convergence of Richardson's
linear solver that stems from the coercivity property \eqref{CoercivityProperty}, in particular
for the case of non-HPD impedance operators. It does converge for Choice 1, but slowly
(relative residual threshold is reached $\mrm{res}(n)<10^{-6}$ after $n=2220$ iterations) which is why we truncated the
convergence history in this case. In addition, this plot exhibits a case where
a non-HPD impedance (Choice 1) is outperformed by an HPD impedance (Choice 2) which is itself
outperformed by another non-HPD impedance (Choice 3 with $\theta = \pi/10$). 

\quad\\
  To conclude, we consider Choice 3 of impedance with values of $\theta$ varying in an interval
around $0$. Figure \ref{RequiredIter} represents the minimum value $n_\theta$
required to obtain $\mrm{res}(n_\theta)<10^{-6}$. We are looking for a value of
$\theta$ that minimizes $n_\theta$. This minimum is reached approximately for $\theta = 0.13$.
It is \textit{not} located at $\theta = 0$ which suggests that considering an imaginary
part for the impedance can be beneficial.


\end{document}